\documentclass[12pt,longbibliography]{article}
\usepackage[utf8]{inputenc}
\usepackage[T1]{fontenc}
\usepackage[a4paper, margin=2.7cm]{geometry}
\usepackage{amsmath, amsthm, amsfonts, amssymb}
\usepackage{mathrsfs}
\usepackage[colorlinks=true,linkcolor=blue,citecolor=blue,urlcolor=blue,breaklinks]{hyperref}
\usepackage{bbm} 
\usepackage{graphicx}
\graphicspath{ {./img/} }
\usepackage{caption}
\usepackage{breakurl}
\usepackage[sort]{natbib}
\usepackage{xcolor}
\usepackage{url}
\def\R{\mathbb{R}}

\def\Z{\mathbb{Z}}

\def\1{\mathbbm{1}}

\renewcommand{\div}{\operatorname{div}}

\def\d{\,\mathrm{d}}

\def \ddt{\frac{\mathrm{d}}{\mathrm{d}t}}

\def\p{\partial}

\def\:{\colon}

\theoremstyle{definition}

\theoremstyle{remark}

\theoremstyle{example}

\numberwithin{equation}{section}
\def\thetitle{Discrete minimizers of the interaction energy in
  collective behavior: a brief numerical and analytic review}
\def\theauthor{José~A.~Cañizo \and Alejandro Ramos-Lora}

\hypersetup{pdftitle={\thetitle}}
\hypersetup{pdfauthor={\theauthor}}
\title{\thetitle}
\author{\theauthor}
\date{April 2024}
\AtEndDocument{\vspace*{12cm}\bigskip{\footnotesize
		\noindent
		\textsc{J.~A.~Cañizo. Departamento de Matemática Aplicada \& IMAG,
			Universidad de Granada.}
		\textit{E-mail address}: \texttt{\href{mailto:canizo@ugr.es}{canizo@ugr.es}}
		\par
		\addvspace{\medskipamount}
		\noindent
		\textsc{A.~Ramos-Lora. Departamento de Matemática Aplicada \& IMAG,
			Universidad de Granada.}
		\textit{E-mail address}: \texttt{\href{mailto:ramoslora@ugr.es}{ramoslora@ugr.es}}
}}

\begin{document}
	
\maketitle

\begin{abstract}
  We consider minimizers of the $N$-particle interaction potential
  energy and briefly review numerical methods used to calculate
  them. We consider simple pair potentials which are repulsive at
  short distances and attractive at long distances, focusing on
  examples which are sums of two powers. The range of powers we look
  at includes the well-known case of the Lennard-Jones potential, but
  we are also interested in less singular potentials which are
  relevant in collective behavior models.  We report on results using
  the software GMIN developed by Wales and collaborators for problems
  in chemistry. For all cases, this algorithm gives good candidates
  for the minimizers for relatively low values of the particle number
  $N$. This is well-known for potentials similar to Lennard-Jones, but
  not for the range which is of interest in collective
  behavior. Standard minimization procedures have been used in the
  literature in this range, but they are likely to yield stationary
  states which are not minimizers. We illustrate numerically some
  properties of the minimizers in 2D, such as lattice structure, Wulff
  shapes, and the continuous large-$N$ limit for locally integrable
  (that is, less singular) potentials.
\end{abstract}

\tableofcontents

\newpage
\section{Introduction}
\label{sec:intro}

Finding configurations with minimum potential energy for a set of
particles is a common mathematical problem found across several
fields. In its simplest version we consider a potential
$V \: \R^d \to \R \cup \{+\infty\}$ and a number $N \geq 2$ of
particles. The problem then consists in finding the configurations
$X = (x_1, \dots, x_N) \in \R^{Nd}$ such that the energy
\begin{equation}
  \label{eq:EN}
  E_N(X) = \sum_{i=1}^N \sum_{\substack{j=1\\j\neq i}}^N V(x_i - x_j)
\end{equation}
is as small as possible. These configurations are called
\emph{minimizers}\footnote{In some contexts they are called
  \emph{global minimizers}, but in the present paper ``minimizer''
  means ``global minimizer'' unless we explicitly mention local
  minimizers.} or \emph{ground states}. In a broader sense, the
problem also involves understanding the structure of these minimizers:
their shape, the distribution of their points, and their asymptotic
properties for large $N$. The energy \eqref{eq:EN} is clearly
invariant by a translation of the points $(x_1,\dots, x_N)$, and
consequently any translation of a minimizer is also a minimizer. One
usually removes this invariance by looking for minimizers with some
fixed property, most often by fixing the mean position
$\sum_{i=1}^N x_i$ to be $0$. If $V$ is radially symmetric, then $E_N$
is also invariant by rotations of the points $(x_1,\dots, x_N$, and by
their inversion $(-x_1, \dots, -x_N)$ (and hence by any rigid motion
of these points). In all cases, minimizers are the same if $V$ is
substituted by $V+C$ for some constant $C$, so some normalization on
$V$ is usually assumed (often, either
$\lim_{|x| \to +\infty} V(x) = 0$ or $V(0) = 0$, but not necessarily).

In this paper we give a short review of numerical methods to calculate
these minimizers, and report on some results in dimension $d=2$ using
the publicly available GMIN software developed by David J. Wales (see
\citet{WalesGMIN}). We carry out numerical simulations in 2D for
potentials of the form
\begin{equation}
  \label{eq:V-powers}
  V(x) = \frac{|x|^a}{a} - \frac{|x|^b}{b},
\end{equation}
where $a > b$. Due to our choice of coefficients, the term with the
$a$ power always represents an attractive potential, and the term with
the $b$ power is repulsive. We will also give a short review of the
available rigorous results to date.

In the previous statement of the problem, the particles are all
identical: they interact with others in an identical way through the
pair potential $V$. More complicated versions of the problem may allow
for particles of different types, interacting in different ways with
one another. One of the better known contexts where this problem is
relevant is chemical physics, where one often wants to find the
configuration of minimum possible energy for a given molecule. This is
important for example in the study of protein folding, where the
natural structure of a protein may be related to its minimum-energy
configuration. In the context of molecular structure, the
Lennard-Jones interaction potential
\begin{equation}
  \label{eq:Lennard-Jones}
  V(x) = \frac{A}{|x|^{12}} - \frac{B}{|x|^6},
\end{equation}
with $A, B > 0$ given constants, is often taken as a test case. The
study of its minimizers (in dimension $d=3$) has received a lot of
attention, and is a notable example of a hard global optimization
problem.

In addition to these problems in physics, this paper is motivated by
several collective behavior models where these minimizers are
relevant. Among them, perhaps the simplest is the so-called
\emph{aggregation equation},
\begin{equation}
  \label{eq:aggregation}
  \p_t u = \div (u (\nabla u * V) ),
\end{equation}
which is a partial differential equation for $u = u(t,x)$, which
represents the density of a certain population which interacts through
the potential $V$. Equilibria of \eqref{eq:aggregation} (if any) must
satisfy $\nabla u * V = 0$ on the support of $u$, which is a property
shared by the minimizers of the \emph{continuous interaction energy}
\begin{equation}
  \label{eq:energy}
  E(u) := \int_{\R^d} \int_{\R^d} u(x) u(y) V(x-y) \d x \d y
\end{equation}
in the set of probability measures (or positive measures with a fixed
total mass). Hence probability minimizers of the energy
\eqref{eq:energy} are particular examples of stationary solutions to
\eqref{eq:aggregation}. To avoid repetition, whenever we speak of
minimizers of $E$ it is understood that we always refer to minimizers
in the set of probability measures. Finding minimizers of $E$ is a
``continuous version'' of the problem of minimizing $E_N$. The link
between these problems is interesting, especially if we consider that
numerical studies of the continuous problem must inevitably carry out
some sort of discretization of the equation, which often leads to some
version of the discrete minimization problem of $E_N$. Hence a proper
understanding of the large-$N$ behavior of minimizers of $E_N$ is also
useful to justify numerical methods for calculating minimizers of the
continuous energy $E$.

There are also many individual-based models whose
``organized states'' share similarities with minimizers of $E_N$. As
an example we mention a model from \citet{DOrsogna2006}: We consider
$N$ individuals with positions and velocities given by $(x_i, v_i)$
for $i = 1, \dots, N$, and the system of ordinary differential
equations
\begin{align*}
  \ddt x_i &= v_i,
  \\
  \ddt v_i &= (\alpha - \beta |v_i|^2) v_i
             - \sum_{\substack{j=1\\j\neq i}}^N \nabla V(x_i-x_j).
\end{align*}
This models a set of individuals which interact through a potential
$V$, but also have a preferred movement speed
$|v| = (\alpha / \beta)^{1/2}$. It is easy to see that if
$X = (y_1, \dots, y_N)$ is a critical point of the energy $E_N$, then
$(x_i(t), v_i(t)) := (y_i + t v, v)$ is a solution to the previous
system of equations. In particular, minimizers of $E_N$ yield coherent
states or \emph{flocking states} for this model. Minimizers (or
critical points) of $E_N$ appear similarly in other models, and this
has motivated interest in types of potentials different from the ones
usually encountered in physics problems.

\paragraph{Existence of minimizers of $E_N$}

For a very general family of potentials it is easy to show that
minimizers must exist. If one assumes that $V$ is lower semicontinuous
and radially strictly increasing whenever $|x| > R$ for some given
$R > 0$, then one may carry out an argument like that in
\citet[Theorem 4.1]{Canizo2018c} to show their existence: it is not
hard to argue that there cannot be ``gaps'' of size larger than $R$
between the particles of a minimizer, which gives a simple estimate on
their diameter. A different argument is given in \citet[Lemma
2.1]{Carrillo2014b} for power-law potentials, and in \citet[Section
1.2]{Blanc2015} for potentials satisfying
$\lim_{|x| \to +\infty} V(x) = 0$.

\section{Numerical methods}
\label{sec:numerical}

In general the energy $E_N$ from \eqref{eq:EN}, as a function
$E_N \: \R^{Nd} \to \R \cup \{+\infty\}$, is not convex, and the
global optimization problem of finding its minima cannot use the
comparatively simpler methods of convex optimization. As mentioned
before, there's the problem of invariance by translations, easily
avoided by fixing the mean position of particles. The problem of
invariance by rotations is harder to avoid, and hence algorithms
usually find one of any possible rotations of the minimizer. It can
then be hard to decide whether two runs of the same algorithm have
given a rotation of the same minimizer, or a fundamentally different
minimizer. However, the deeper problem is that there are usually a
huge number of stationary states for the energy $E_N$; that is, points
$(x_1, \dots, x_N)$ at which
\begin{equation}
  \sum_{\substack{j=1\\j\neq i}}^N \nabla V(x_i - x_j)
  =
  0
  \qquad \text{for all $i = 1, \dots, N$.}
\end{equation}
Many of these are local minima, and runs of standard local
minimization algorithms will most often fall into one of these states,
and not get to a minimizer. Hence simple descent algorithms such as
variants of conjugate gradient descent, or Newton-type methods like
LBFGS \citep{liu1989limited}, will almost surely not find minimizers even for small
$N$. Global minimization algorithms must somehow avoid these ``traps''
by more clever search techniques. We will give a short overview of
numerical methods used for this kind of problems, and then describe
the algorithm we use, implemented in the program GMIN by D.~J.~Wales
\citep{WalesGMIN}.

Most numerical methods found in the literature have been tested with
the Lennard-Jones potential \eqref{eq:Lennard-Jones}, and most often
in 3D, which is a challenging problem of physical interest and has
been the subject of many publications in chemistry journals. A
database with all currently available minima of $E_N$ for several
potentials in 3D (and molecular minimization problems of other kinds)
is available online \citep{CCDonline}, among them Lennard-Jones and
Morse potentials. This database contains references of the first
publication for each minimizer, and the reader may also check the
unpublished paper \citet{Zhang2010} for a good list of historical
references.

\begin{figure}
  \centering
  \includegraphics[width=14cm]{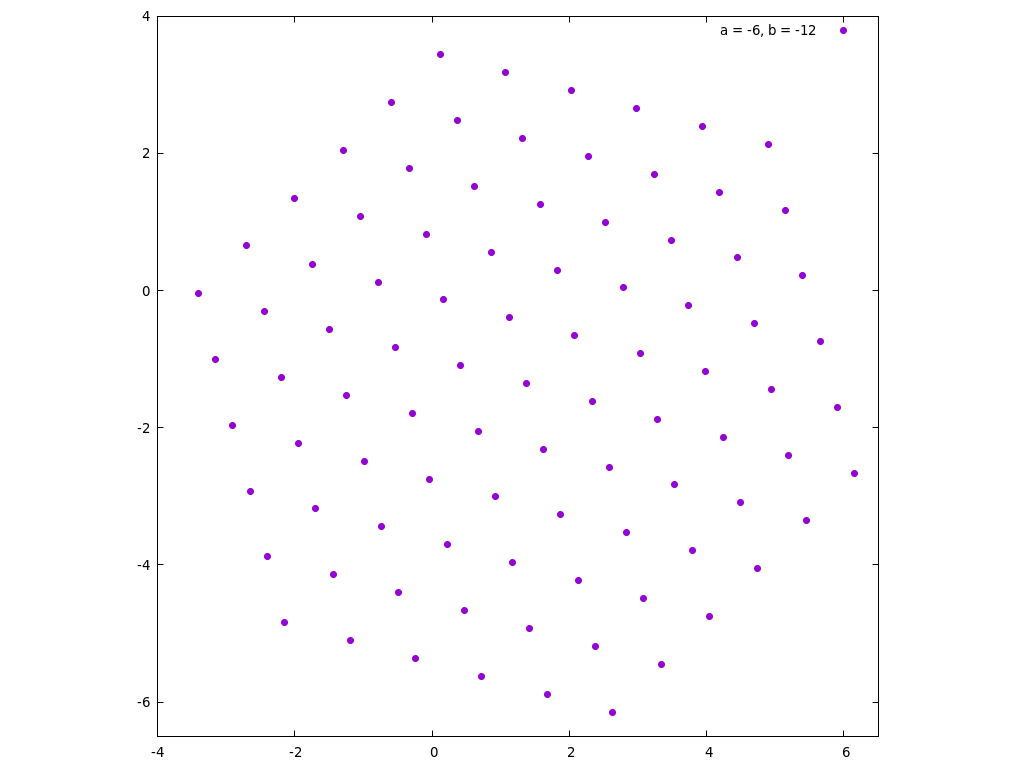}
  \caption{Minimizer for the Lennard-Jones potential with 91
    particles.}
  \label{fig:Lennard-Jones}
\end{figure}

Early methods such as those in
\citet{Hoare1979,HOARE1972,Hoare1971a,Hoare1971} are \emph{biased}
methods, in the sense that they work by ``guessing'' that minimizers
have some kind of regular crystal structure. Initial guesses for a
minimizer are obtained by placing particles at the sites of an
appropriate lattice, and then running local minimization algorithms to
optimize the positions. Then one carries out some kind of random
perturbation of these positions in order to find ``nearby'' local
minima which may have a lower energy. In addition to guessing initial
candidates by using the sites of a lattice, one may also start by
finding minimizers for small $N$, and add particles to build
candidates for larger $N$. These methods were quite successful in
finding many candidates for minimizers which are still standing today,
but there are several ``hard'' values of $N$ for which better (lower
energy) configurations were found later \citep{Northby1987,
  Wales1997}. An important step of many of these methods, used for the
``random perturbation'' we mentioned, is \emph{simulated annealing}
\citep{Kirkpatrick1983}, a global optimization method which has an
appealing physical motivation for chemistry problems. Roughly
speaking, in simulated annealing a random perturbation $X'$ of a given
configuration $X$ is considered, and then:
\begin{enumerate}
\item It is accepted if $X'$ has lower energy than $X$,
  or
\item It is accepted with a probability proportional to
  $\exp(-\frac{E_N(X') - E_N(X)}{T})$ if $X'$ has higher energy than
  $X$.
\end{enumerate}
Here $T$ is a parameter which has an analogy to temperature, and
governs how probable it is to jump ``far'' from the current state. In
the literature on this, performing this algorithm any number of times
for a fixed $T$ is sometimes referred to as a \textsc{Monte Carlo
  simulation}, and each step as a \emph{Monte Carlo step}. Simulated
annealing methods then work by running many steps of this random
perturbation, and gradually reducing $T$ until departures from the
current state become very unlikely. Methods focused on this technique
have had some success \citep{Wille1987a,Wille1987}, and are also part
of more recent methods.

Another family of methods is that of evolutionary algorithms, which
have also been quite successful in finding minima for Lennard-Jones
potentials, even for the values of $N$ considered harder
\citep{Daven1996,Niesse1996}. Genetic algorithms are a particular kind
of evolutionary algorithms, though the naming is not completely
precise. These methods are based on an analogy with natural selection,
and are also an example of biased methods which make informed guesses
on the structure of minimizers. We refer to the more recent paper by
\citet{Dittner2017} for a list of references on the use of this kind
of algorithms for Lennard-Jones energy minimization.

The paper \citet{ben2003numerical} contains a good presentation of
related minimization problems, and describes a simple algorithm for
minimizing $E_N$ for the Lennard-Jones potential, which is tested
in 2D. It is also a biased method in the sense that particles are
initially placed at the sites of a triangular lattice, and later
suitable symmetrization methods are used.

\begin{figure}
  \centering
  \includegraphics[width=15cm]{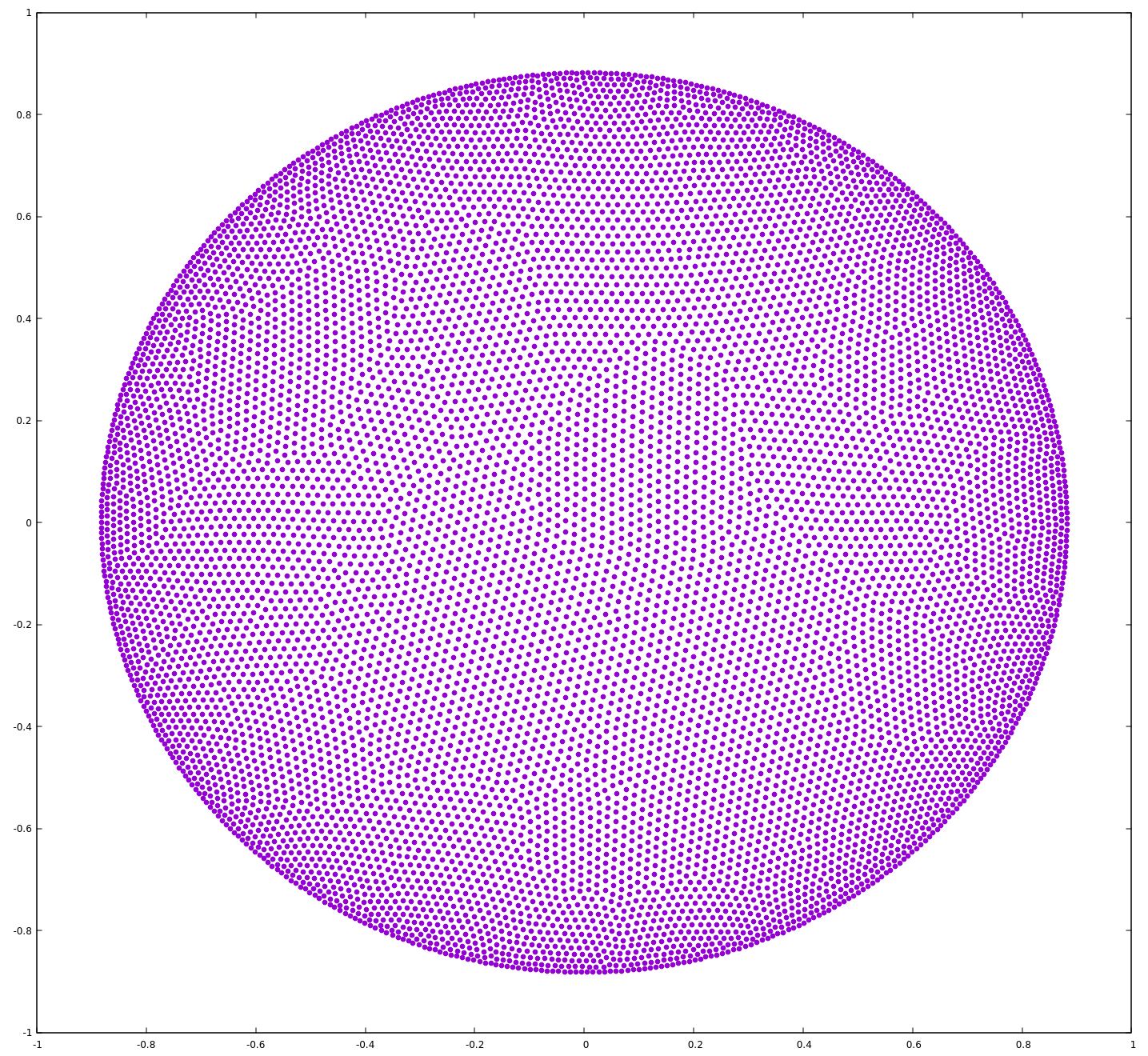}
  \caption{Our best numerical candidate for a minimizer of the
    potential \eqref{eq:V-powers} with $a=2$, $b=0.5$ and $10^4$
    particles. (We thank David J.~Wales for providing this result)
  }
  \label{fig:a2b05-10000}
\end{figure}

\begin{figure}[h!]
  \centering
  \includegraphics[width=14cm]{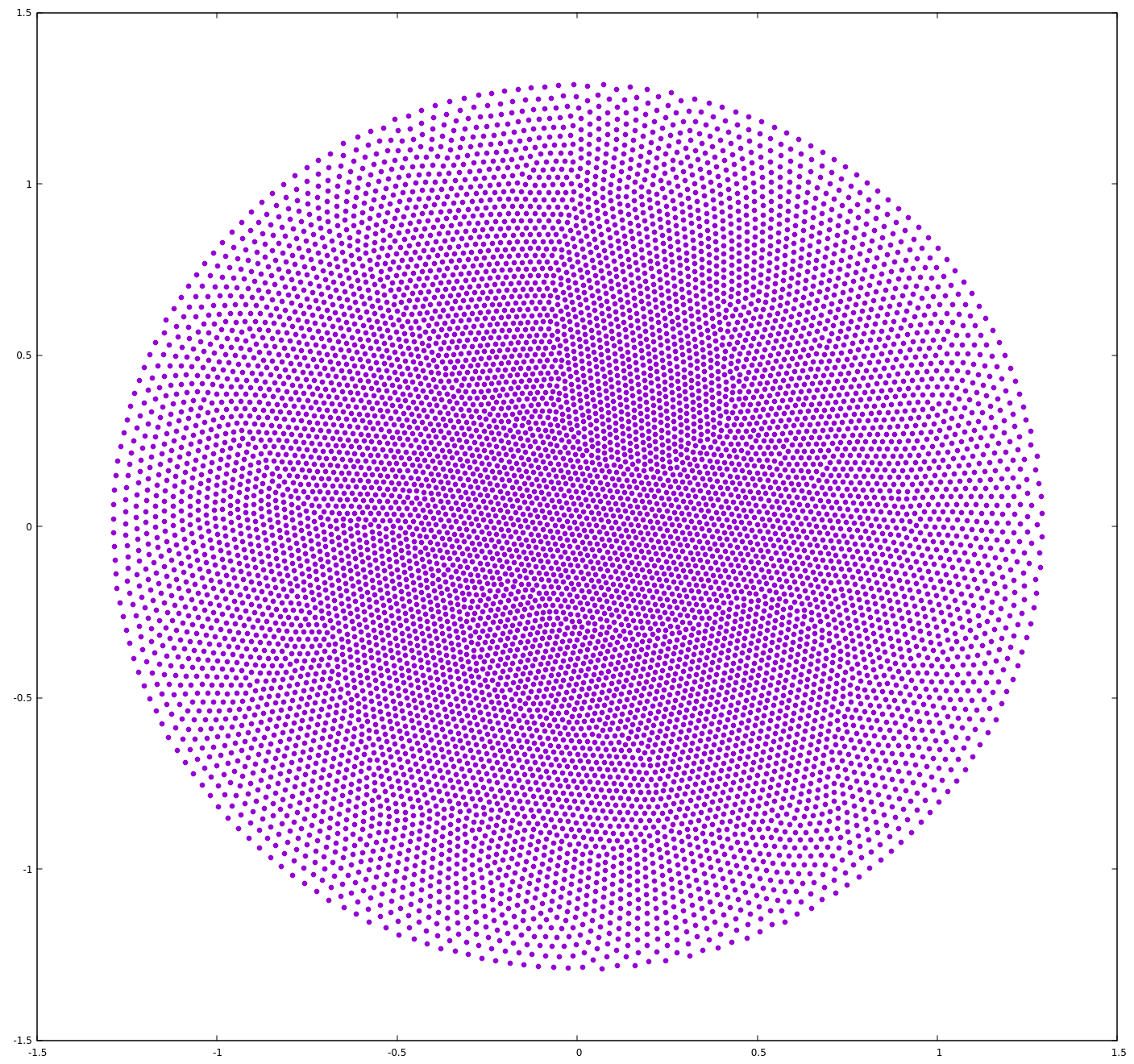}
  \captionsetup{singlelinecheck=off}
  \caption[]{Numerical candidate for minimizer of the potential
    \eqref{eq:V-powers} with $a=2$, $b=-1$ and $10^4$ particles. (We
    thank David J.~Wales for providing this result.) This corresponds to
    the ``strongly repulsive'' regime considered in
    \citet{Carrillo_2016}. In this regime it is conjectured that a
    sequence of minimizers $(X_N)_{N \geq 2}$ must converge to a
    measure $\nu$ which must have a smooth density with respect to
    Lebesgue measure. In fact, the most likely candidate measure is
    given in \citet[eq. (38)]{A_Carrillo_2017} and is of the form
    \begin{equation*}
      \nu(x) = C (R^2 - |x|^2)_+^{1 - \frac{b+d}{2}},
    \end{equation*}
    with $(\cdot)_+$ denoting the positive part, with suitable
    constants $C, R > 0$ which ensure in particular that the total
    mass is $1$. In the case corresponding to this figure we have
    $d=2$, $b=-1$, so the density of $\nu$ vanishes at the
    boundary. By contrast, in the case of Figure \ref{fig:a2b05-10000}
    the exponent $1 - (b+d)/2$ becomes negative and the density of
    $\nu$ should diverge at the boundary.}
  \label{fig:a2b-1}
\end{figure}

As discussed, these biased methods have been tried mainly for
potentials of Lennard-Jones type, motivated by problems in chemistry
and crystallography. Since our aim is to use a good numerical
algorithm for a wider range of potentials whose minimizers may behave
very differently, an adaptation of these methods does not seem
straightforward, especially for the cases which have a continuous
limit such as those discussed in Section \ref{sec:continuous}, where
particles do not arrange as large-scale crystals. Instead, we use a
technique known as \emph{basin-hopping}, first reported for the
Lennard-Jones potential in \citet{Wales1997}, and before used by
\citet{Li1987} for a different problem in chemistry. This technique
has been in particular quite successful for Lennard-Jones type
potentials. It is an \emph{unbiased} algorithm in the sense that it
does not assume any particular structure of the minimizers, and it
seems a more appropriate algorithm for the range of potentials we are
interested in. It works as follows: first, the energy $E_N(X)$ is
substituted by an energy $\widetilde{E}_N(X)$, obtained by performing
a suitable descent algorithm starting from the configuration $X$ (the
Polak-Ribiere variant of the conjugate gradient algorithm is mentioned
in \citet{Wales1997}, although it currently uses a customized LBFGS
method). That is: if we denote by $Y = M(X)$ the result of running a
descent algorithm from the starting configuration $X$, then
\begin{equation*}
  \widetilde{E}_N(X) := E_N( M(X) ).
\end{equation*}
Hence the new \emph{energy landscape} $\widetilde{E}_N$ is a piecewise
constant function, which is constant in a region where our chosen
descent algorithm stops at the same local minimum $Y$, and whose value
in this region is $E_N(Y)$. This new landscape can be studied
theoretically, taking $M(X)$ to be the point to which the descent
algorithm converges, and it can be approximated numerically in
practice by considering the point at which the numerical descent
stops. Now several steps are followed:
\begin{enumerate}
\item A Monte Carlo exploration of the landscape $\widetilde{E}_N$ is
  carried out with a fixed temperature, chosen so that the acceptance
  rate of new points is about $1/2$ (as explained in our brief earlier
  discussion of simulated annealing).
  
\item During this exploration, we also observe the \emph{individual
    energy} (or just energy) of each particle $i$, defined by
  $$E_{N,i} = \sum_{j \neq i} V(x_i - x_j).$$
  We consider the particle $i_{max}$ which has a maximum energy among
  all $N$ particles, and the one $i_{min}$ which has minimum
  energy. If this maximum energy is larger than a certain multiple of
  the lowest energy, then an \emph{angular move} is performed for the
  particle $i_{max}$: this particle is removed from the cluster, and
  placed at a random position close to the boundary of the cluster.
\end{enumerate}
In addition to this, a further Monte Carlo exploration is carried
out with the following special starting positions:
\begin{enumerate}
\item The best configuration found by the algorithm with $N-1$
  particles, adding one particle at the best position found by a fixed
  number of angular moves.
  
\item The best configuration found by the algorithm with $N+1$
  particles, removing the particle with highest energy.
\end{enumerate}
We have applied this algorithm for different potentials, otherwise
with no modifications. Some of the results are explored in the next
sections, together with a short review of the expected behavior of
minimizers. We would like to emphasize that while obtaining local
minimizers is straightforward in almost all cases, obtaining good
candidates for a \emph{minimizer} is a hard numerical problem. As an
example, the hexagonal-shaped minimizers for very singular potentials
such as Lennard-Jones (see Fig \ref{fig:Lennard-Jones}) require some
kind of global search algorithm as described above. Or consider for
example a minimizer for the locally integrable potential
\eqref{eq:V-powers} with $a=2$ and $b=0.5$ (Figure
\ref{fig:a2b05-10000}). Asymptotically as $N \to +\infty$, this
density is known to converge to the density of a minimizer of the
continuous energy \eqref{eq:energy} (see Section
\ref{sec:continuous}), but the way this happens shows several
interesting phenomena. One sees that locally particles still seem to
be arranged in triangular patterns; but since the overall shape must
be round (since the support of the minimizer is strongly suspected to
be a ball), these local triangular arrangements must undergo some kind
of dislocation patterns in order to fit the global shape. In Figure
\ref{fig:a2b05-10000} one can clearly appreciate dislocation curves
along which triangular patterns with different orientations
meet. These patterns are quite sensitive to the minimization algorithm
used, and we do not know whether the plot in Figure
\ref{fig:a2b05-10000} actually represents a minimizer. But an
understanding of these patterns is certainly an interesting problem,
and even its numerical study requires powerful algorithms which can
handle a large number of particles and yield good approximations to
the minimizer. For particle numbers over a few thousand it is likely
that this is still out of reach, even numerically. In Figures
\ref{fig:a2b05-10000} and \ref{fig:a2b-1} we show minimizer candidates
for 10,000 particles, which have required large computational power
and several weeks of computer time.

\section{Crystallization phenomena}
\label{sec:crystal}

The wider problem of understanding the reason why matter is often
arranged in crystalline patterns is known as the problem of
crystallization. The reviews by \citet{Ventevogel1978,Radin1987} are a
good introduction to this, giving a motivation and related
mathematical problems. It happens that minimizers of $E_N$ have a
close relationship to the observed arrangement of atoms or molecules
at low temperature, and hence an aspect of the crystallization problem
which is particularly simple to state is: \emph{prove rigorously that
  (or whether) minimizers of $E_N$ approach a periodic arrangement as
  $N \to +\infty$}. For some potentials $V$ which have a singularity
stronger than $|x|^{-d}$ at $x=0$, minimizers of $E_N$ numerically
seem to do exactly this, but proofs are surprisingly difficult even in
dimension $1$! The review by \citet{Blanc2015} is a nice source of
references on this problem, and we would like to give a short summary
with our current understanding.

Let us first state the problem a bit more precisely, essentially
following \citet{Blanc2015}. In the context of crystallization, we
usually consider potentials $V = V(x)$ which are positive at $x=0$ (or
have a singularity where the potential diverges to $+\infty$ as
$x \to 0$), have a strictly negative minimum value, and satisfy
\begin{equation}
  \label{eq:V-lim0}
  \lim_{|x| \to +\infty} V(x) = 0.
\end{equation}
Further conditions are needed for minimizers to show any kind of
crystallization behavior; see more on this below. The canonical
example of the potentials which exhibit this kind of phenomena is the
Lennard-Jones potential from equation \eqref{eq:Lennard-Jones}.

For a particle distribution $X = (x_1, \dots, x_N) \in \R^{Nd}$ we
consider the associated \emph{empirical measure} given by
\begin{equation}
  \label{eq:muX}
  \mu_X := \sum_{i=1}^N \delta_{x_i}.
\end{equation}
Notice that we are \emph{not} dividing by $N$, so this measure has
total mass equal to $N$. In many cases we expect this measure to
approach some periodic arrangement as $N \to +\infty$, and in order to
define this precisely we recall that a \emph{lattice} in $\R^d$ is a
subgroup of $\R^d$ of the form
\begin{equation}
  \label{eq:lattice}
  L = \{ k_1 v_1 + k_2 v_2 + \dots + k_d v_d \mid k_1, \dots, k_d \in \Z \},
\end{equation}
where $(v_1, \dots, v_d)$ form a vector space basis of
$\R^d$. Lattices are a simple form of periodic structures, and are a
basic concept in crystallography, where they are often referred to as
\emph{Bravais lattices}. The reason for this special naming is that in
crystallography the name \emph{lattice} is often broadened to include
periodic structures which are not lattices in the sense of the above
definition; notable examples are the honeycomb arrangement in 2D, or
the hexagonal close packed arrangement in 3D. Here we will always use
the word \emph{lattice} in the mathematical sense stated in
\eqref{eq:lattice}. Consistently with \eqref{eq:muX}, for a lattice
$L$ we denote by $\mu_L$ the measure
\begin{equation}
  \label{eq:muL}
  \mu_L := \sum_{p \in L} \delta_{p}.
\end{equation}
Notice that $\mu_L$ has infinite mass (is an unbounded measure), but
it is finite in compact sets. Here is a very strong formulation of the
problem of crystallization:
\begin{description}
\item[\textbf{Q1 (lattice crystallization)}] Does there exist a
  lattice $L \subseteq \R^d$ and a sequence of minimizers
  $(X_N)_{N \geq 2}$ of $E_N$ such that (possibly up to a subsequence)
  \begin{equation*}
    \mu_{X_N} \overset{*}{\rightharpoonup} \mu_L
    \qquad \text{as $N \to +\infty$?}
  \end{equation*}
  The convergence here is the weak-$*$ convergence of measures
  (convergence when tested against continuous, compactly supported
  functions). Or even stronger: does this take place for \emph{all}
  sequences of minimizers, up to rigid movements of these minimizers?
\end{description}
Crystalline structures may well be periodic structures which are not
lattices. A very natural formalization of a periodic structure is a
measure $\mu$ on $\R^d$ such that $\mu(\cdot + p) = \mu(\cdot)$ for
all $p$ in some lattice $L$ (that is, a measure which is invariant by
translations of vector $p$, for all $p$ in a given lattice). If this
happens, we say that $\mu$ is \emph{invariant by the lattice
  $L$}. Following this, here is a weaker form of the above problem:
\begin{description}
\item[\textbf{Q2 (periodic crystallization)}] Does there exist a
  lattice $L \subseteq \R^d$ and a sequence of minimizers
  $(X_N)_{N \geq 2}$ of $E_N$ such that (possibly up to a subsequence)
  \begin{equation*}
    \mu_{X_N} \overset{*}{\rightharpoonup} \mu
    \qquad \text{as $N \to +\infty$,}
  \end{equation*}
  where $\mu$ is a measure invariant by the lattice $L$? Or even
  stronger: does this take place for \emph{all} sequences of
  minimizers, up to rigid movements of these minimizers?
\end{description}
If these questions are too hard, one may consider the large-$N$
behavior of some macroscopic quantities like the energy per
particle. Given a lattice $L$ in $\R^d$ and a potential $V$ we define
its \emph{energy per particle} by
\begin{equation}
  \label{eq:energy-per-particle-lattice}
  e(L) := \sum_{p \in L\setminus\{0\}} V(p),
\end{equation}
when this sum converges absolutely. Also, given a configuration $X = (x_1,
\dots, x_N) \in \R^{Nd}$ we define its energy per particle by
\begin{equation*}
  e_N(X) := \frac{1}{N}E_N(X).
\end{equation*}
Here is a related, usually weaker, statement of the problems
\textbf{Q1} and \textbf{Q2}:
\begin{description}
\item[\textbf{Q3 (energy crystallization)}] Does there exist a
  lattice $L \subseteq \R^d$ and a sequence of minimizers
  $(X_N)_{N \geq 2}$ of $E_N$ such that (possibly up to a subsequence)
  \begin{equation*}
    e_N(X_N) \to e(L)
    \qquad \text{as $N \to +\infty$?}
  \end{equation*}
  Or even stronger: does this take place for \emph{all} sequences of
  minimizers?
\end{description}
In a related problem, one can try to give good estimates for the
convergence of $e_N(X_N)$ to $e(L)$, for example giving bounds on its
rate of convergence. The fact that $e_N(X_N)$ has \emph{some} limit is
known for a wide range of potentials.

In the context of crystallization phenomena one always assumes
\eqref{eq:V-lim0}. As a consequence, a sequence of configurations with
$N$ particles such that the interparticle distances all diverge to
$+\infty$ has energy as close to $0$ as we wish. Hence for a minimizer
$X_N$ it must happen that $E_N(X_N) \leq 0$, so if \textbf{Q3} holds
one must have $e(L) \leq 0$. In any case, if \textbf{Q3} holds for the
full sequence of minimizers (and not up to a subsequence), in
particular it implies that for some $C > 0$
\begin{equation}
  \label{eq:H-stability}
  E_N(X_N) \geq -C N, \quad \text{for all}\; N \geq 2.
\end{equation}
Potentials $V$ satisfying \eqref{eq:H-stability}
are called \emph{H-stable}, or sometimes just \emph{stable} in the
literature; see for example \citet[Section 3.2]{Ruelle1969}. Often
\eqref{eq:H-stability}, though seemingly weaker, actually implies that
$e_N(X_N)$ does have a limit as $N \to +\infty$. Since the terminology
comes from statistical mechanics and thermodynamics, interest is
placed on potentials for which thermodynamical quantities (such as the
energy-per-particle) have a finite limit, which justifies the word
\emph{stable}. Minimizers for potentials which are \emph{not} stable
may still have different scaling limits which are also interesting,
though their study has different motivations and is more recent; see
Section \ref{sec:continuous} below.

Finally, there is the related question of the \emph{formation of a
  macroscopic object}. The \emph{diameter} of a minimizer is defined as
the maximum distance between two of its points.  If we rescale
minimizers in order to fix their diameter, does the resulting measure
converge to something? If we look at the regime where one expects that
particles approach a periodic distribution, we would expect the
diameter of a minimizer to be roughly proportional to $N^{1/d}$, and
hence we pose the following question:
\begin{description}
\item[\textbf{MO (formation of a macroscopic object)}] Does there exist
  a nontrivial measure $\mu$ on $\R^d$ and a sequence of minimizers
  $(X_N)_{N \geq 2}$ of $E_N$ such that (possibly up to a subsequence)
  \begin{equation}
    \label{eq:macro-convergence}
    \frac{1}{N} \mu_{X_N / N^{1/d}} \overset{*}{\rightharpoonup} \mu
    \qquad \text{as $N \to +\infty$?}
  \end{equation}
  Or even stronger: does this take place for \emph{all} sequences of
  minimizers, up to rigid movements of these minimizers?
\end{description}
The empirical measure $\mu_{X_N / N^{1/d}}$ corresponds to the scaled
minimizer
$$\frac{1}{N^{1/d}}X_N := \left(
\frac{x_1}{N^{1/d}}, \dots, \frac{x_N}{N^{1/d}} \right),$$ and we
divide the measure by $N$ to ensure a total mass equal to $1$. This
question concerns behavior at a different scale from questions
\textbf{Q1}--\textbf{Q3}. An answer to \textbf{MO} seems to not give
much information on any of \textbf{Q1}--\textbf{Q3}, and
reciprocally.

We have focused on the above questions in order to narrow the
discussion to a representative case, since the field is too broad for
a short review. However, there are many related aspects of these
problems. Some of them involve infinite distributions of
particles. Among these we mention the following:
\begin{enumerate}
\item (\emph{Lattice minimization}) Among all lattices in $\R^d$, which ones
  minimize the energy per particle
  \eqref{eq:energy-per-particle-lattice}?
  
\item (\emph{Infinite particle distributions}) Given a countably
  infinite set of points in $\R^d$, it is possible to define their
  \emph{mean energy per particle} in a reasonable way, so that it is
  consistent with the energy per particle of a lattice given in
  \eqref{eq:lattice}. Among all infinite particle distributions, which
  ones minimize the energy per particle?
\end{enumerate}
We do not elaborate on different mathematical interpretations of the
problem, namely: minimization problems with fixed density; modified
models with a background energy; statistical mechanical or quantum
mechanical models. We refer the reader to the reviews mentioned above
for more on these matters \citep{Radin1987,Blanc2015}. For
minimization on \emph{compact} submanifolds of $\R^d$ we refer to
\citet{Hardin2005, Borodachov2019}.

Minimizers for a specific case of potential $V$ have links to sphere
packing problems. If we consider the following \emph{sticky hard
  sphere} potential
\begin{equation}
  \label{eq:sticky}
  V(x) :=
  \begin{cases}
    +\infty \qquad & \text{if $|x| < 1$,}
    \\
    -1 \qquad & \text{if $|x| = 1$,}
    \\
    0 \qquad & \text{if $|x| > 1$,}
  \end{cases}
\end{equation}
then minimizers cannot have any two particles at a distance closer
than $1$ (otherwise their energy would be $+\infty$), and there is a
reward for spheres which are touching. Hence a minimizer of this
potential is just a configuration of $N$ hard spheres in
$d$-dimensional space with the maximum total number of contact points
among them. In dimension $2$ this problem was solved by
\citet{Heitmann1980}; see also \citet{Luca2018} for a more recent
proof. But the problem is far from understood in dimensions $3$ and
higher, with open questions even for low values of $N$ (above $11$). A
popular review of the problem can be found in \citet{Hayes2012}, and
see \citet{Meng2010,Hoy2012} for some recent results. Attempts have
been made to relate the structure of these minimizers to those for the
Lennard-Jones potential \citep{Trombach}, but a precise link remains
unclear. It is also tempting to make a connection of the sticky hard
sphere problem with that of \emph{sphere packing}, which is the
problem of finding a distribution of identical hard spheres in $\R^d$
with maximum possible density. This link has often been remarked, but
as far as we know there are no rigorous results backing this up. There
have been some surprising recent results in sphere packing by
\citet{Viazovska2017,Cohn2017}, who give specific lattices in
dimensions $8$ and $24$ which can be proved to be the densest possible
sphere distributions in these dimensions. The proof relies on the
existence of specially symmetric lattices in these dimensions, which
are natural candidates for this densest packing.

Now that we have made an effort to state rigorously the results we are
interested in, we can try to summarize the available rigorous results:

\paragraph{One-dimensional results}

In dimension $1$ there were several early results by
\citet{Ventevogel1978a,Ventevogel1979a,Ventevogel1979}. These results
concern mainly infinite particle distributions (as mentioned in points
1,2 above) but specifically \citet{Ventevogel1978a} contains a proof
of the strong version of \textbf{Q3} for a wide class of power-law
potentials of the form \eqref{eq:V-powers} with $b < 1$ and a certain
restriction on $a$: the energy per particle of any sequence of
minimizers converges to the energy per particle of a given lattice;
this lattice is the unique one with minimum energy per particle. This
is a rather complete answer to \textbf{Q3} for some power-law
potentials. However, as far as we know there are no satisfying general
conditions on $V$ which ensure a similar answer.

In contrast, problem \textbf{Q1} is not so well understood even in
dimension $1$. The best available result seems to be the one by
\citet{Gardner1979}, which does show the strong version of \textbf{Q1}
for some potentials including the Lennard-Jones one. The conditions
for this to hold are more restrictive, and it seems a hard problem to
generalize them to a wide class of potentials.

In dimension $1$ the minimization problem for the sticky hard sphere
potential \eqref{eq:sticky} is trivial, and minimizers are
distributions of $N$ equally spaced particles at interparticle
distance $1$.

\paragraph{Two-dimensional results}

For the sticky hard sphere potential \eqref{eq:sticky} in dimension
$2$ results are already interesting. A first proof that minimizers
must be placed at the sites of a triangular lattice was given by
\citet{Heitmann1980}, giving an answer to \textbf{Q1} (and hence
\textbf{Q2}) and \textbf{Q3}. For a specific piecewise linear
potential not too far from \eqref{eq:sticky}, a positive answer to
\textbf{Q3} was also given by \citet{Radin1981}. Then \textbf{MO} was
answered by \citet*{AuYeung2011} for a family of short-range
potentials with $V$ compactly supported and close enough to the sticky
sphere potential \eqref{eq:sticky}: the convergence
\eqref{eq:macro-convergence} takes place up to a subsequence, and the
measure $\mu$ is $\frac{2}{\sqrt{3}} \1_E$, where $\1_E$ is the
characteristic of a finite-perimeter set $E$ of area
${\sqrt{3}}/{2}$. In \citet{Theil2006}, \textbf{Q3} is answered
positively for a certain family of potentials which allow similar
growth conditions as the Lennard-Jones one. In addition, it is proved
that with periodic boundary conditions, the triangular lattice is a
minimizer. Recently in \cite{Betermin2023} a computer-assisted method
is used to show that the triangular lattice minimizes energy among
lattices for the Lennard-Jones potential in 2D (hence giving an answer
to the related lattice minimization problem we mentioned earlier).

Adding a three-body potential (and hence for an energy not
strictly of the form \eqref{eq:EN}), and still with periodic boundary
conditions, \citet{Weinan2008} followed the ideas in \citet{Theil2006}
to show that minimizers must have particles placed at the sites of a
\emph{hexagonal} periodic arrangement (the vertices of a honeycomb,
which are not a lattice in the sense we are using in this
paper). Hence in this case, a version of \textbf{Q2} holds, but not
\textbf{Q1}. Problem \textbf{Q3} has been solved by
\citet{B_termin_2021} for some specific piecewise-constant potentials
which give rise to a square lattice.

\paragraph{Three-dimensional results}

From the point of view of rigorous results, the situation in three
dimensions is almost completely open. As far as we know there are no
results regarding any of \textbf{Q1}, \textbf{Q2}, \textbf{Q3} or
\textbf{MO}, even for the sticky hard sphere potential
\eqref{eq:sticky}. The problem of optimal packing, however, was solved
by \citep{Hales2005} with a non-conventional proof which requires
several computer checks. One of the difficulties of the 3D problem is
that there are many non-equivalent sphere distributions in $\R^3$
which maximize density, including non-periodic structures.

\paragraph{Higher dimensional results}

In higher dimensions there have been important breakthroughs on
optimal sphere packing, for which the optimal lattices have been
recently identified in dimension $8$ \citep{Viazovska2017} and $24$
\citep{Cohn2017}. As a consequence of these ideas, in \citet{Cohn2022}
it was proved that these optimal lattices are also \emph{universally
  optimal}: they minimize the interaction energy among infinite
configurations of points with fixed density $1$, for any potentials
$V$ which are completely monotone functions of the squared
distance. It is conjectured that the same holds in dimension $2$ for
the usual triangular lattice, but this is not proved as far as we
know. In \citet{Petrache} it is shown that this universal optimality
also implies a corresponding result for Riesz and Coulomb interaction
energies, which require an appropriate renormalized definition for
infinite particle distributions. Hence this latter result is known in
dimensions $8$ and $24$, and is a conjecture in dimension $2$.

All of these results concern infinite particle configurations. As far
as we know, no results are known concerning any of \textbf{Q1},
\textbf{Q2}, \textbf{Q3} or \textbf{MO} in dimensions $d \geq 3$.

\section{Continuous limits}
\label{sec:continuous}

Consider a sequence $(X_N)_{N \geq 2}$ of minimizers of $E_N$, and a
locally integrable potential $V$ (in the case of potentials of the
form \eqref{eq:V-powers} this holds if and only if $b > -d$). Then it
may happen that
\begin{equation*}
  \nu_{X_N} := \frac{1}{N} \mu_{X_N}
  =
  \frac{1}{N} \sum_{i=1}^N \delta_{x_i^N}
  \overset{*}{\rightharpoonup} \nu
  \qquad \text{as $N \to +\infty$},
\end{equation*}
for some measure $\nu$, where $X_N = (x_1^N, \dots, x_N^N)$. This type
of limit is of \emph{mean-field} type due to the factor $1/N$, and
notice that there is no rescaling as there was in
\eqref{eq:macro-convergence}. In this case the discrete interaction
energy $E_N$ is an approximation of the continuous one $E$,
\begin{multline}
  \label{eq:EN-to-E}
  \frac{1}{N^2} E_N(X)
  =
  \frac{1}{N^2}
  \sum_{i=1}^N \sum_{\substack{j=1\\j\neq i}}^N V(x_i - x_j)
  =
  \int_{\R^d} \int_{\R^d} V(x-y) \d \nu_{X_N}(x) \d \nu_{X_N}(y)
  \\
  \approx
  \int_{\R^d} \int_{\R^d} V(x-y) \d \nu(x) \d \nu(y)
  = E(\nu),
\end{multline}
when this can be rigorously justified. Then we expect the measure
$\nu$ to be a minimizer of the continuous interaction energy
\eqref{eq:energy}. Proving this is another interesting question which
we will state explicitly:
\begin{description}
\item[MF (Mean-field limit)] Does there exist a nontrivial measure
  $\nu$ in $\R^d$ and a sequence of minimizers $(X_N)_{N \geq 2}$ of
  $E_N$ such that (possibly up to a subsequence)
  \begin{equation*}
    \frac{1}{N} \mu_{X_N} \overset{*}{\rightharpoonup} \nu
    \qquad \text{as $N \to +\infty$?}
  \end{equation*}
  Does this take place for \emph{all} sequences of minimizers, up to
  rigid movements of these minimizers? In this case, is $\nu$ a
  minimizer of the continuous interaction energy \eqref{eq:energy}?
\end{description}
In analogy to \textbf{Q3}, one may also look for good asymptotics of
the energy-per-particle $E_N(X_N)/N$ as $N \to +\infty$. For locally
integrable potentials, if \textbf{MF} holds we expect
\eqref{eq:EN-to-E} to hold and hence we also expect that
\begin{equation*}
  E_N(X_N) \sim N^2 E(\nu).
\end{equation*}
It is interesting to study this approximation independently, and
improve it to more precise approximations. This problem was studied by
\citet{Petrache2017} when the potential $V$ is repulsive and there is
an additional external confining potential. For the particular case
$-d < b < a=2$ in \eqref{eq:V-powers}, and when continuous minimizers
are regular enough, this gives quite precise asymptotics for
$E_N(X_N)$.

The only answer we know to question \textbf{MF} was given in
\citet*{Carrillo2014b} for power-law potentials \eqref{eq:V-powers}
with $b \geq 1$, and in \citet{Canizo2018c} for general potentials
which include the power-laws \eqref{eq:V-powers} with $b > 2-d$. In
this case, the answer to all questions in \textbf{MF} is positive for
\emph{strictly unstable} potentials $V$, that is, potentials $V$ for
which there exists a probability distribution $\rho$ on $\R^d$ such
that
\begin{equation*}
  E(\rho) < \lim_{|x| \to +\infty} V(x)
\end{equation*}
(assuming this limit exists or is $+\infty$). For a wide range of
potentials, this concept is actually shown to imply that the classical
H-stability condition \eqref{eq:H-stability} does \emph{not} hold
\citep{Simione2015,Canizo2018c}. We believe \textbf{MF} should also be
true in the range $-d < b \leq 2-d$, but as far as we know there is no
available proof. For the strategy in \citet{Canizo2018c} two
ingredients are needed:
\begin{enumerate}
\item One needs to show that the sequence of discrete energies $E_N$
  Gamma-converges in a suitable sense to the continuum energy
  $E$.
  
\item One needs to show that the diameter of a minimizer $X_N$ stays
  bounded as $N \to +\infty$. As a consequence, the empirical measures
  $\mu_{X_N} / N$ converge up to a subsequence, in the weak-$*$
  sense, to a certain probability measure $\mu$. By using the previous
  Gamma-convergence result (and possibly some additional uniform
  estimates on other properties of the minimizers) one may then deduce
  that $\mu$ must be a minimizer of $E$.
\end{enumerate}

\medskip In relation to these questions, a description of the
minimizers of $E$ is also an interesting problem. As remarked in the
introduction, these minimizers are of interest for collective behavior
models, and are also quite interesting from the point of view of
calculus of variations. A general result on existence of these
minimizers was given in \citet{Canizo2015b,Simione2015}, but there
remain many questions on their uniqueness, symmetry and regularity
properties. There has been a lot of recent activity regarding these,
and we will try to give here a short account of the literature on
this.

An early discussion on the problem of minimizing $E$ can be found in
\citet[Section X]{Bavaud1991}.  In \citet{Raoul2012}, for very regular
potentials it is proved that stationary states (and in particular
minimizers, if any) must be sums of delta functions. A link between
the repulsive singularity of the potential at $x=0$ and the dimension
of the support of minimizers was given by \citet*{Balague2013},
essentially by noticing that for any minimizer $\rho$ the convolution
$\Delta V * \rho$ must be bounded. A consequence of this is that the
support of a minimizer for $V$ of the form \eqref{eq:V-powers} with
$2-d \leq b < 2$ must have dimension at least $2-b$. This is not
optimal since numerical simulations seem to indicate that, at least
for ``usual'' potentials, the support always has integer dimension,
and some effort has been devoted to understanding this effect. See
however \citet{Carrillo_2022} for an example of fractal behavior for a
specifically tailored potential $V$.

For ``strongly repulsive'' potentials with $-d < b < 2-d$ (but still
locally integrable), regularity of minimizers can be obtained via
obstacle problems \citep*{Carrillo_2016}. When $a=2$ and $-d< b<2-d$,
these minimizers correspond to the ``fractional Barenblatt solutions''
of \citet{Caffarelli2011a, Caffarelli2011}; candidates for minimizers
can be found for $-d < b < a$ with $b \leq 2$, when $a$ or $b$ are
even integers \citep{A_Carrillo_2017}. In dimension $d \geq 2$, when
$a=2$ and $2-d<b<4-d$, minimizers are radially symmetric, unique, and
were obtained in \citet{Carrillo_2022}.

For $b = 2$ and $a \in (2,4)$, the unique minimizer (up to
translations) is the uniform measure on an $(n-1)$-dimensional sphere
of appropriate radius \citep*{Davies2021}. This has recently been
extended by \cite{frank2023minimizers} to the case $d \geq 2$,
$2 \leq a \leq 4$, $b_*(a) \leq b \leq 2$ with $b < a$, for a certain
decreasing function $b_* = b_*(a)$ which is always between $-d+3$ and
$-d+4$.

This is in contrast with minimizers for \emph{mildly repulsive}
potentials (potentials which are $\mathcal{C}^2$ at $0$), for which
minimizers must be supported on a finite number of points under quite
general conditions \citep*{Carrillo_2017}. This result is true even
for \emph{local} minimizers in the topology of the
$\infty$-Wasserstein transport distance. This applies to power-law
potentials of the form \eqref{eq:V-powers} whenever $2 < b < a$. For
$a \geq 4$ and $b \geq 2$ (excluding the case $(a,b) = (4,2)$) the
minimizer is unique up to rigid movements, and it must be a sum of
equal delta functions supported at the vertices of a regular
$d$-simplex \citep*{Davies_2022}. In this same reference, the case
$(a,b) = (4,2)$ is studied. It is a special limit case in which
minimizers are not unique (even up to rigid movements): minimizers may
be supported on a sphere, on a finite number of points, or may be
convex combinations of these possibilities. In the 1D case explicit
minimizers are known for some power-law potentials of the form
\eqref{eq:V-powers} \citep{Frank2022}.

In addition to this, the question of stability or instability of these
minimizers for several dynamical models has been investigated in a few
cases; see for example \citet{Fellner2010,Fellner2010a} (for 1D
stability of sums of Dirac masses for the aggregation equation
\eqref{eq:aggregation}); \citet{Balague2013a, Kolokolnikov2011} (stability of
spherical shells and rings); \citet{Bertozzi2012} (behavior of patch
solutions for $b=2-d$).

We will also mention some results on non-isotropic
potentials: see \citet{Carrillo2021a} and \citet*{Carrillo2023a,
  Carrillo2022}.

There are also related minimization problems which would take us too
far from the topic; we refer the reader to \citet{Blanc2015} and
\citet{Frank_2023} for an introduction to similar problems and many
references on them.

\begin{figure}[h!]
  \centering
  \includegraphics[width=3.5cm]{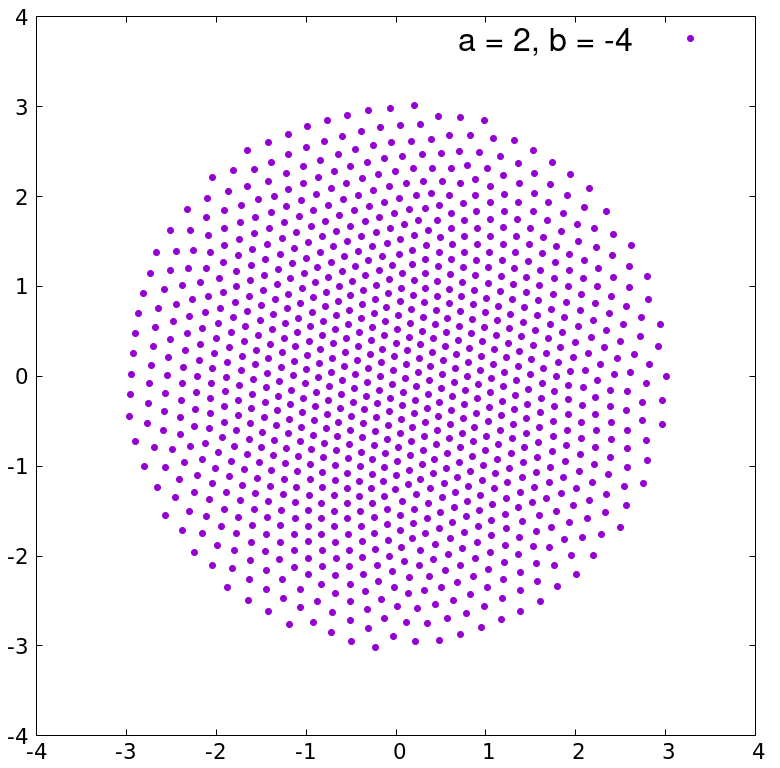}
  \includegraphics[width=3.5cm]{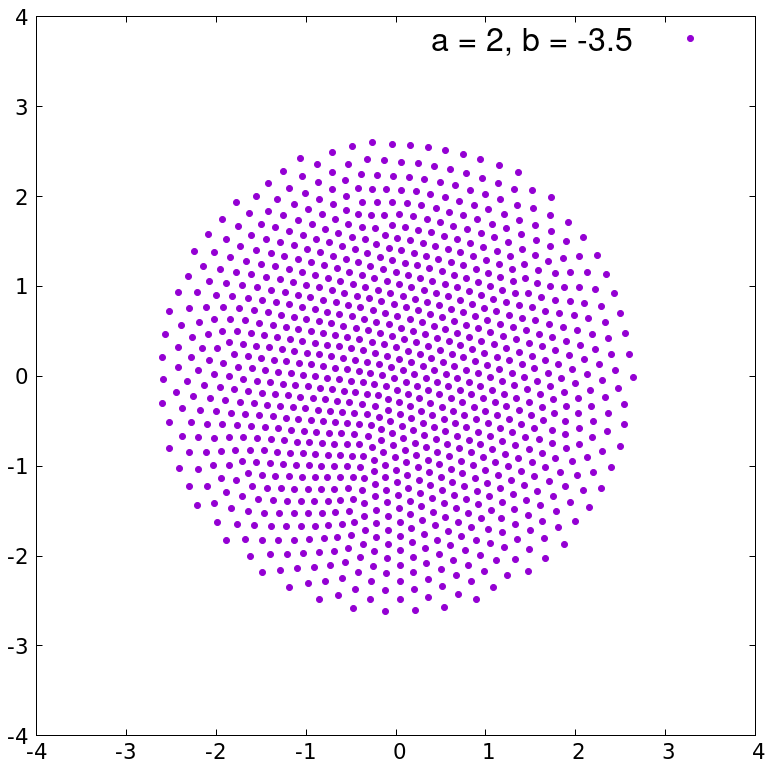}
  \includegraphics[width=3.5cm]{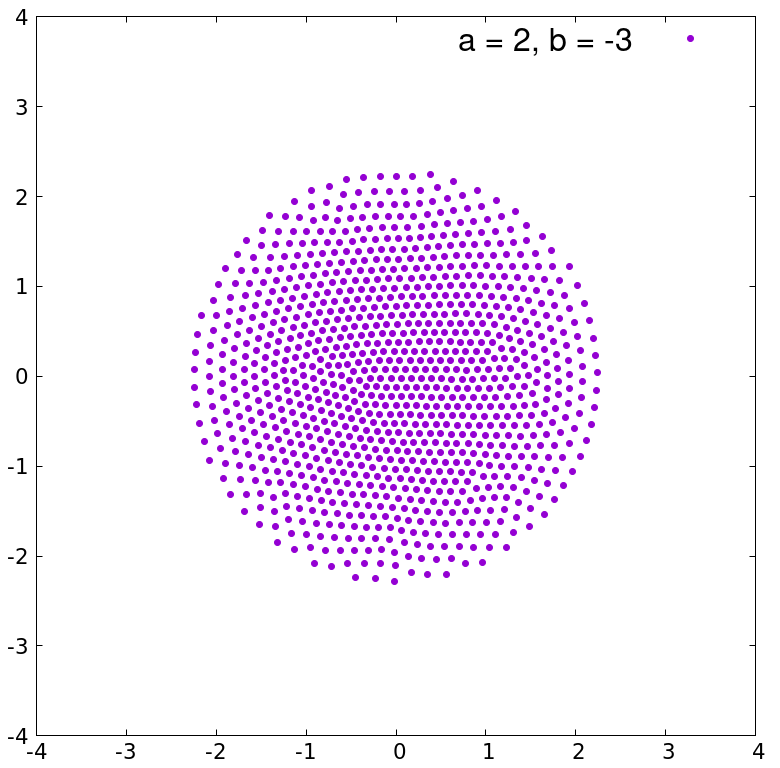}
  \includegraphics[width=3.5cm]{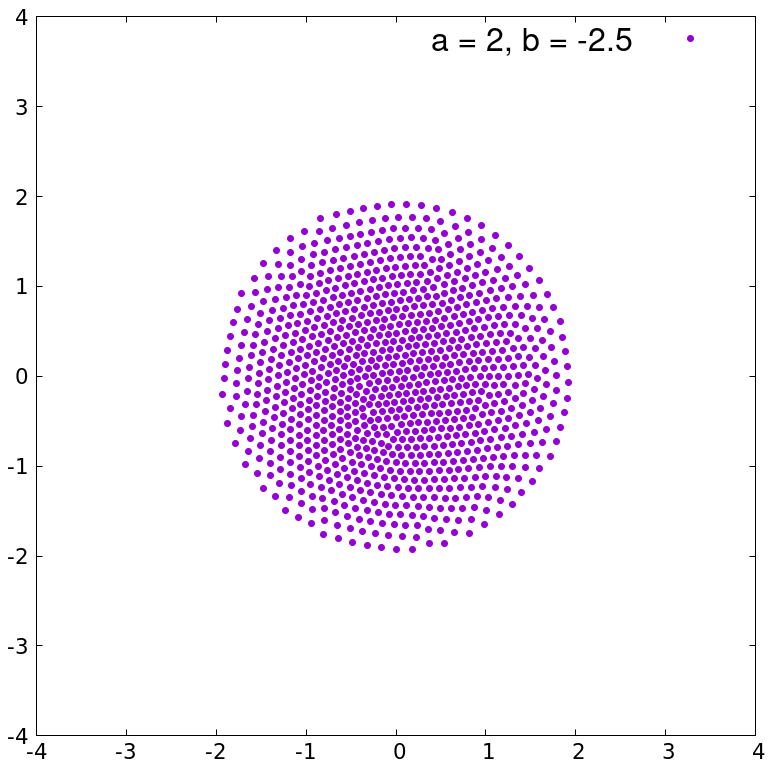}
  \includegraphics[width=3.5cm]{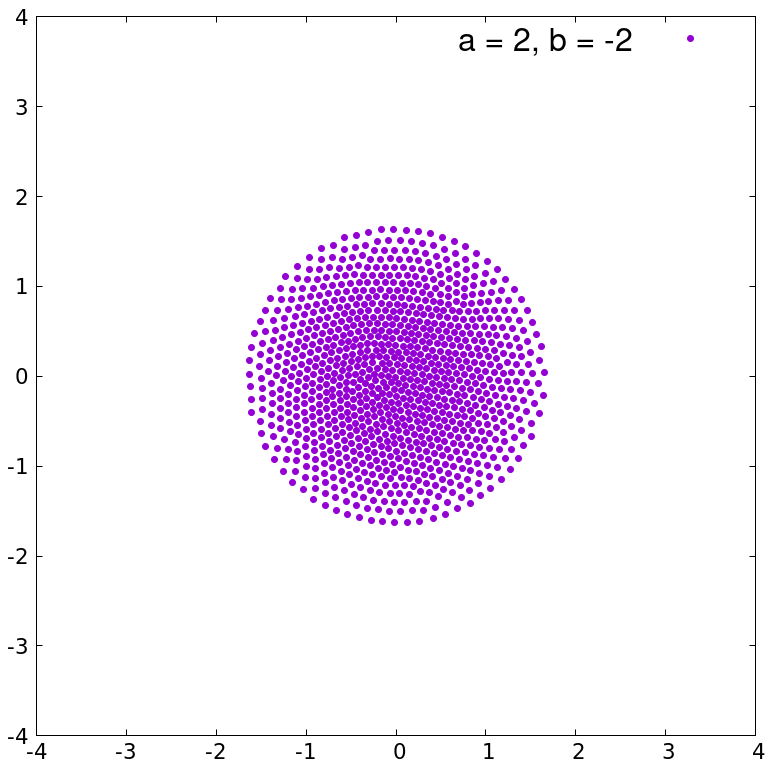}
  \includegraphics[width=3.5cm]{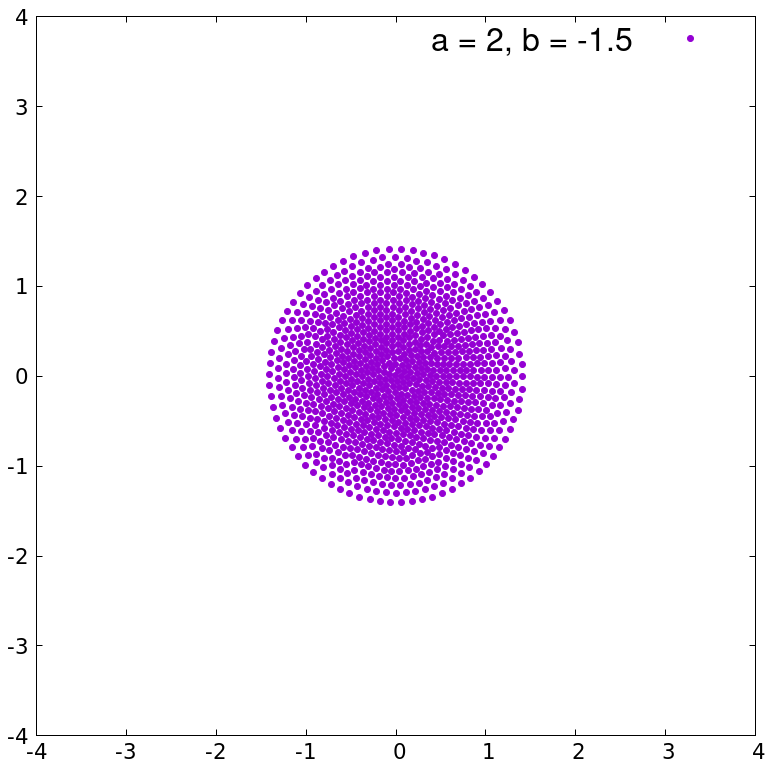}
  \includegraphics[width=3.5cm]{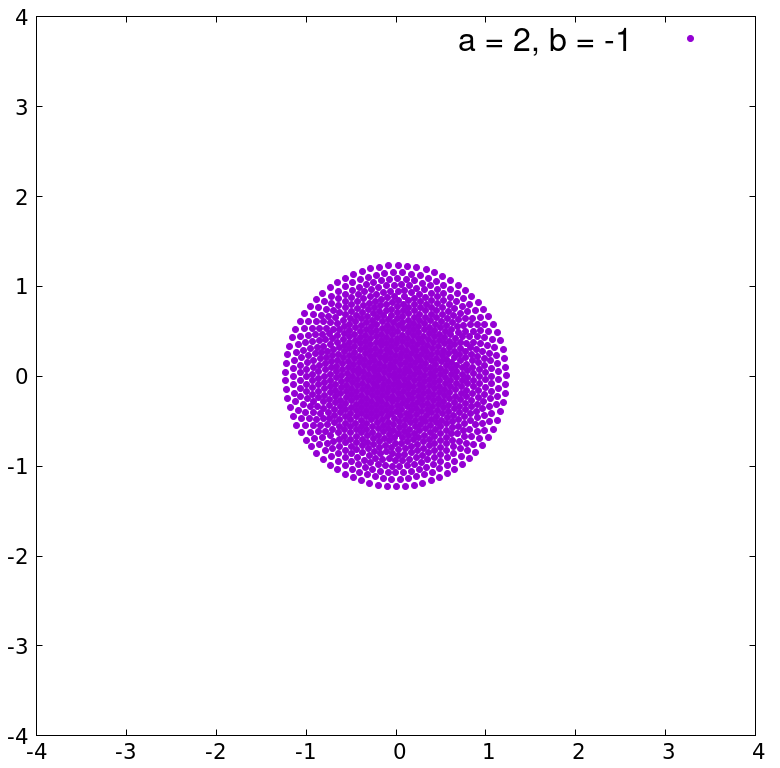}
  \includegraphics[width=3.5cm]{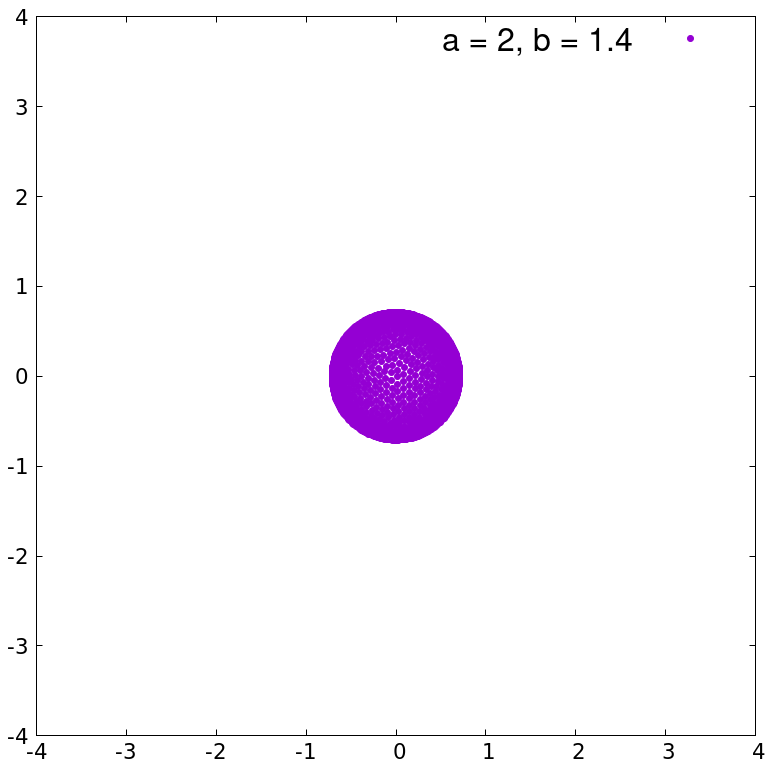} 
  \caption{Minimizer for the power-law potential \eqref{eq:V-powers}
    with $a=2$ and several values of $b$, always with $10^3$
    particles.}
  \label{fig:a2}
\end{figure}

\begin{figure}[h!]
  \centering
   \includegraphics[width=7cm]{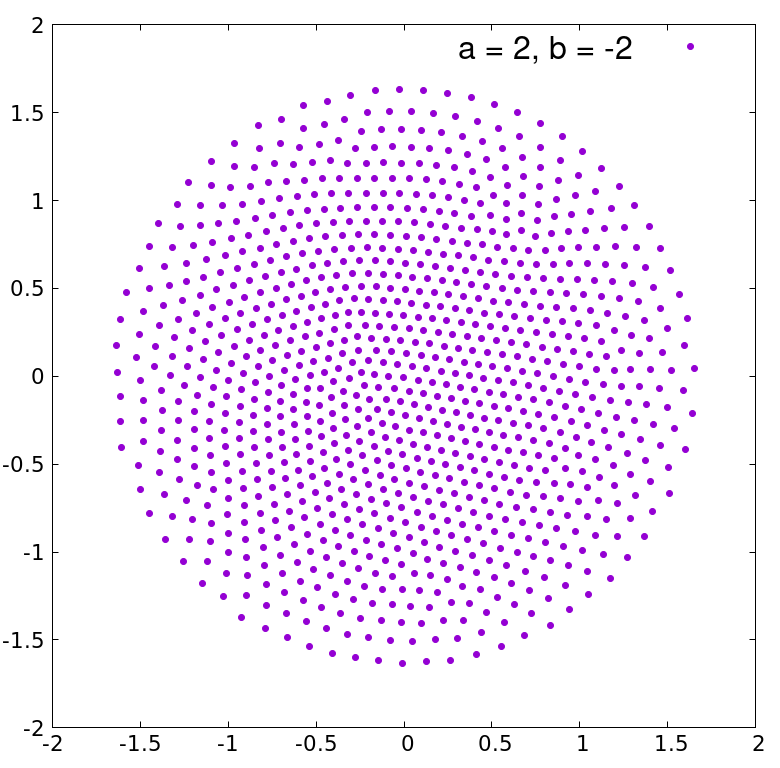}
   \includegraphics[width=7cm]{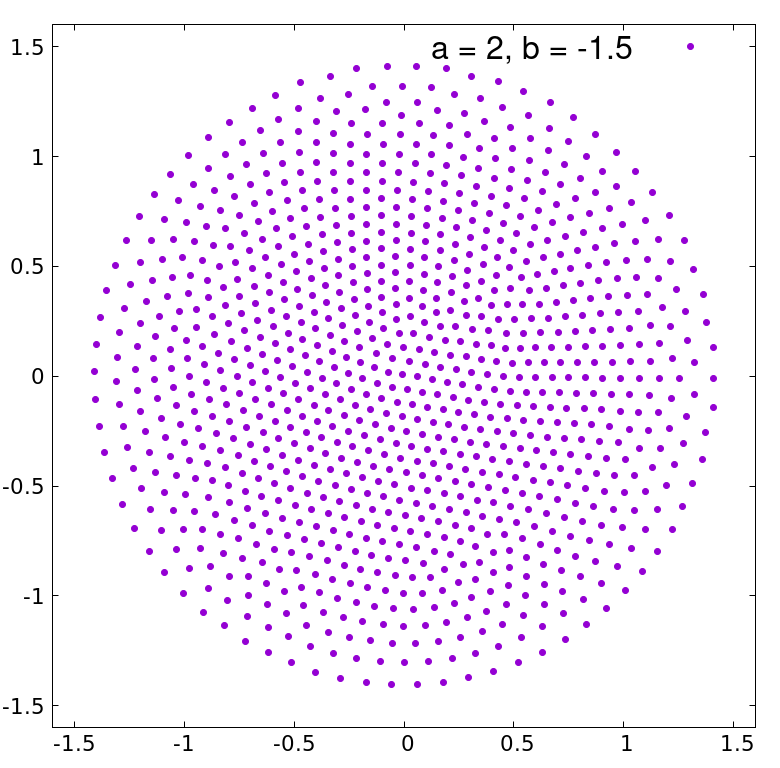}
   \includegraphics[width=7cm]{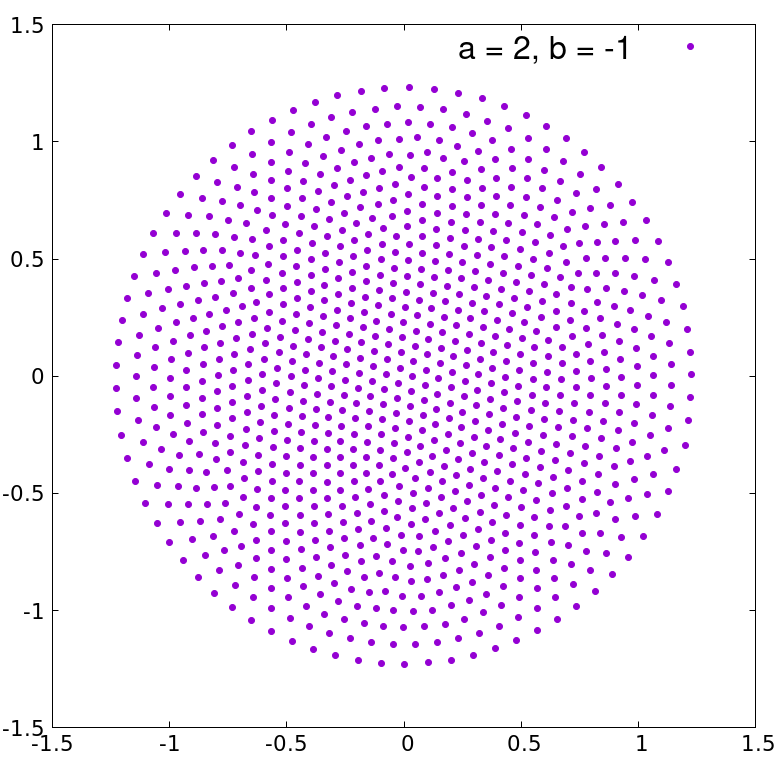}
   \includegraphics[width=7cm]{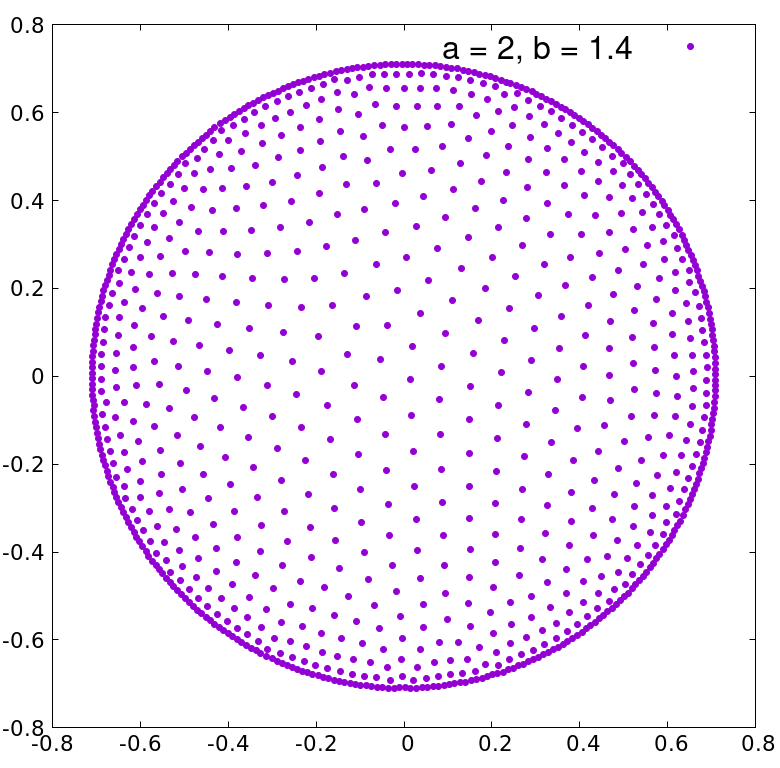}
  \caption{Zoomed-in view of the last 4 minimizers in Figure
    \ref{fig:a2}.}
  \label{fig:a2-zoom}
\end{figure}

We will illustrate the above properties in 2D with a few numerical
results for power-law potentials of the form \eqref{eq:V-powers}. Let
us first fix $a=2$ and look at the appearance of numerical minimizers
for several values of $b$. As the exponent $b$ increases, the
repulsion becomes progressively weaker. It is clearly seen in Figures
\ref{fig:a2} and \ref{fig:a2-zoom} that minimizers become more
concentrated as $b$ increases, since repulsion becomes locally
weaker. Also, they become more concentrated on the boundary of a ball
for higher $b$. For $b < 0$ we are in the regime studied in
\citet{Carrillo_2016}, so we know minimizers of $E$ should be smooth,
but it is not proved that $N$-particle minimizers approach minimizers
of $E$. In contrast, for $b > 0$ we know from \citet{Canizo2018c} that
$N$-particle minimizers must approach minimizers of the continuous
energy $E$, but we do not know the specific shape of these minimizers;
for example, it is not proved in this range that minimizers must be
regular functions on their support, or even $L^P$ functions for some
$p \geq 1$. Reasonable candidates for minimizers of $E$ are given in
\citet{A_Carrillo_2017}, and they are all smooth functions on their
support.

\begin{figure}[h!]
  \centering
   \includegraphics[width=3.5cm]{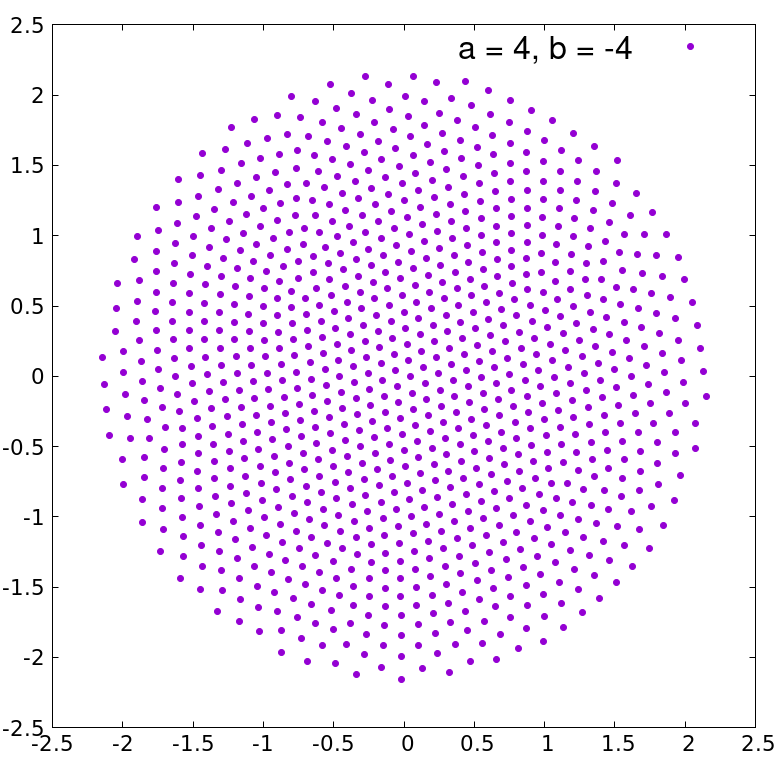}
   \includegraphics[width=3.5cm]{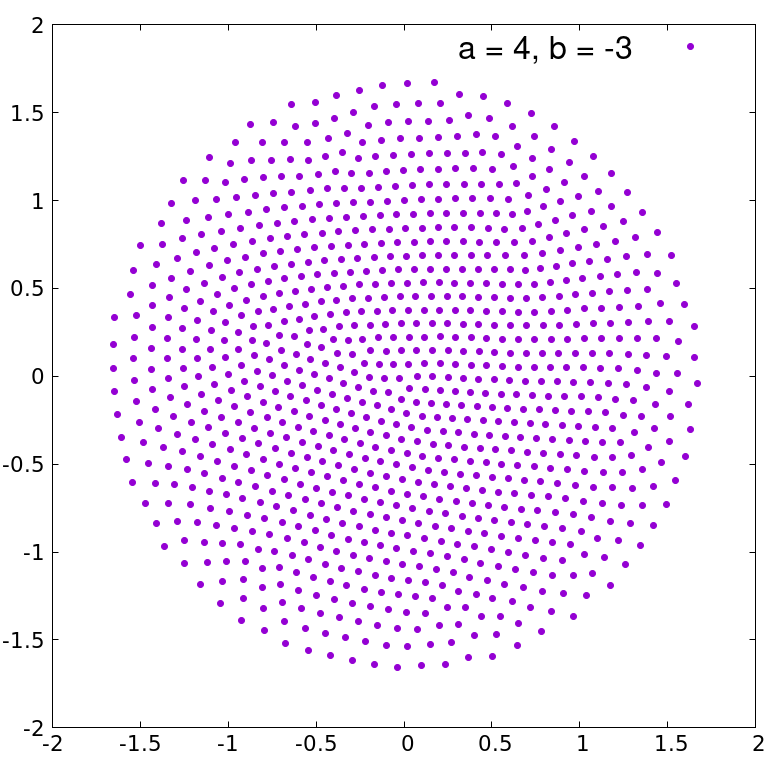}
   \includegraphics[width=3.5cm]{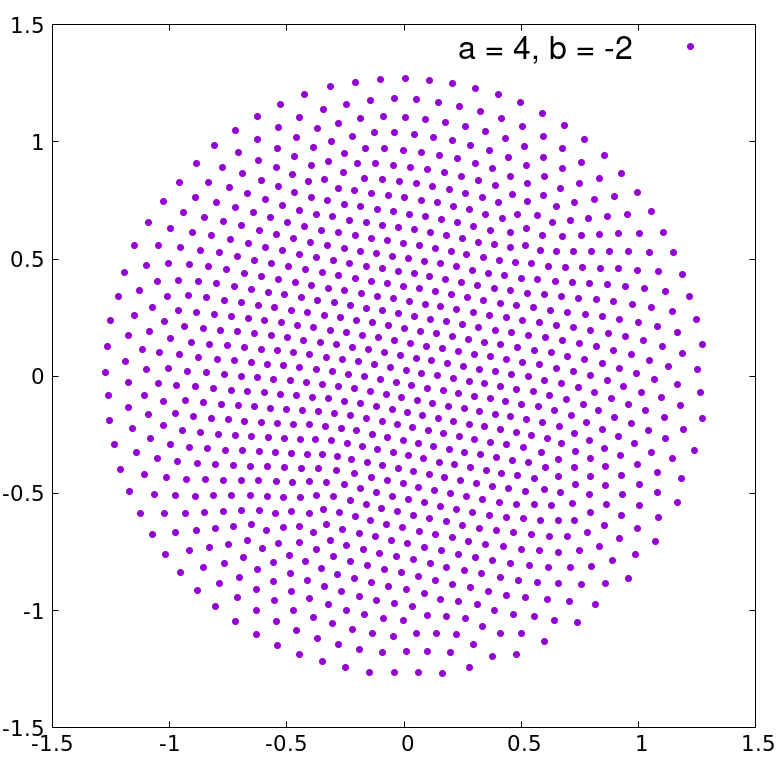}
   \includegraphics[width=3.5cm]{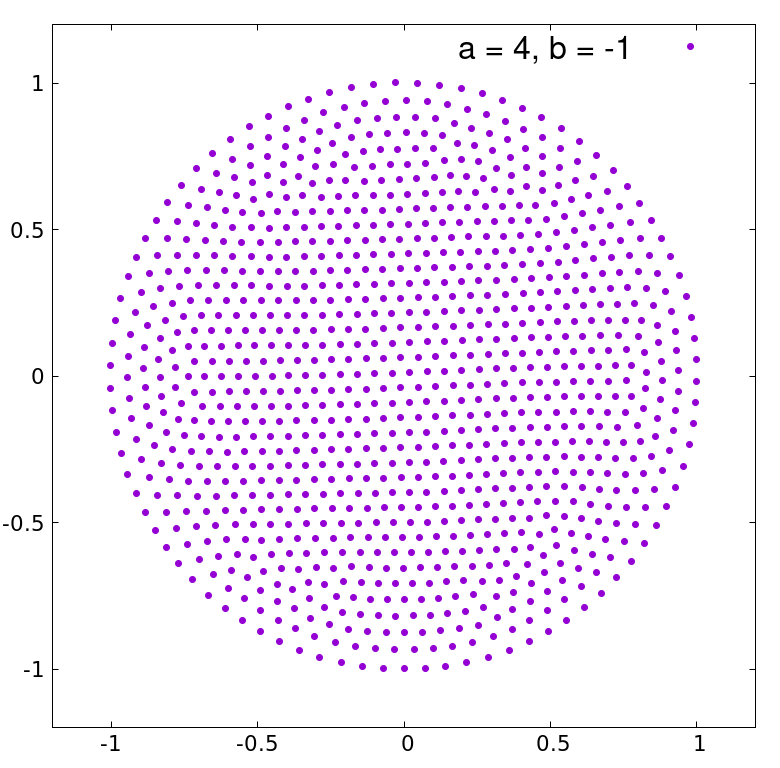}
   \includegraphics[width=3.5cm]{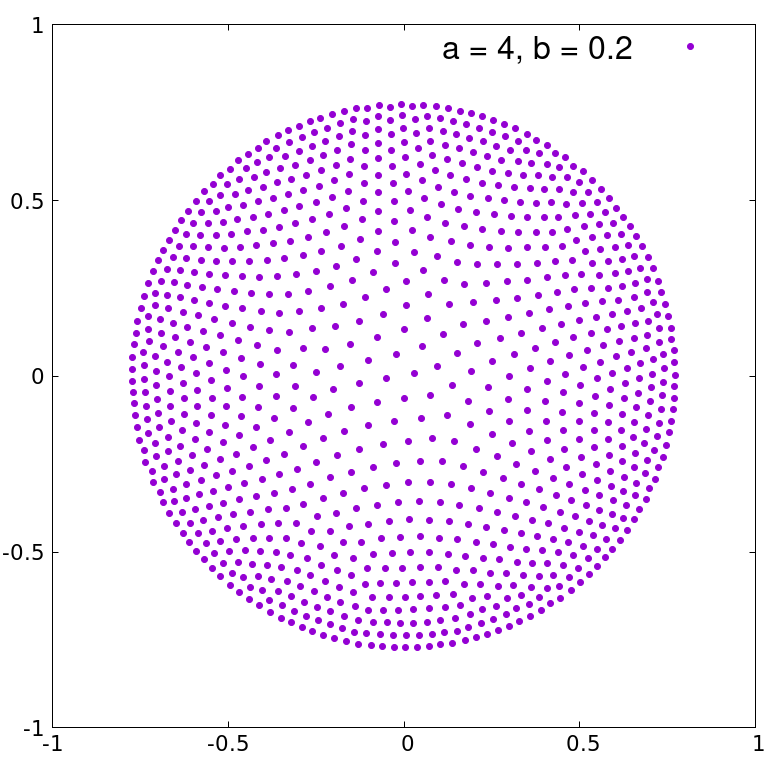}
   \includegraphics[width=3.5cm]{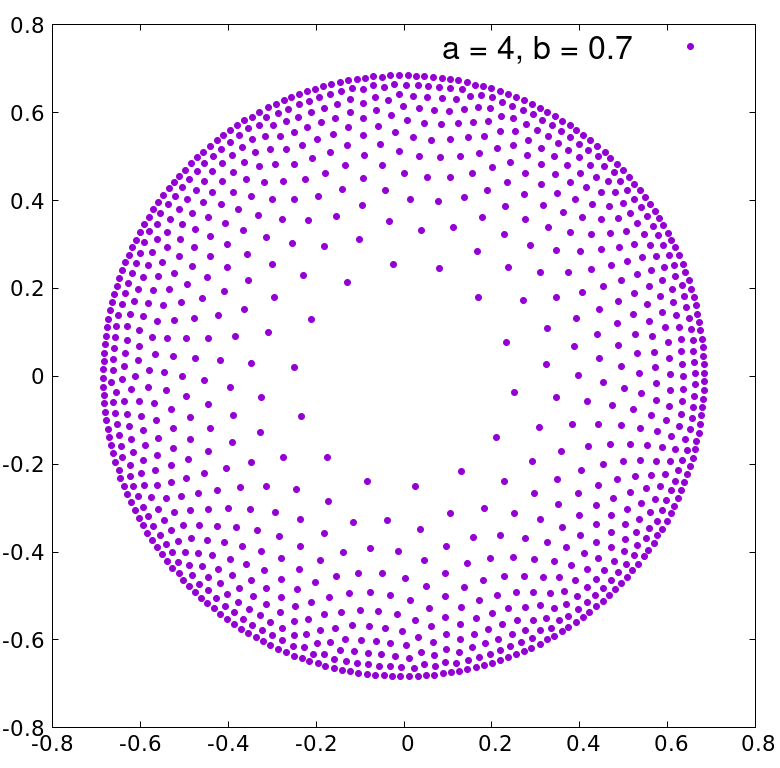}
   \includegraphics[width=3.5cm]{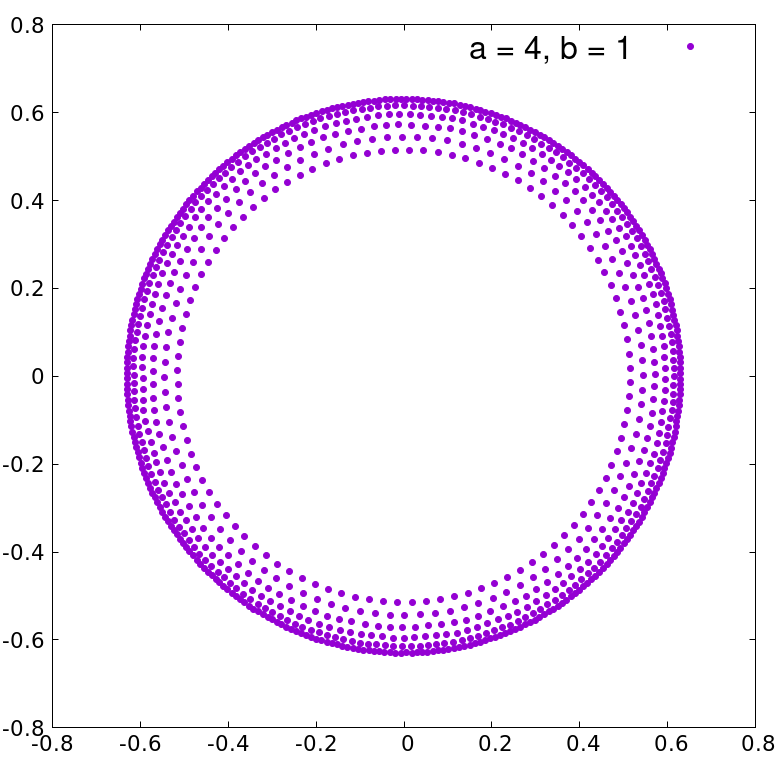}
   \includegraphics[width=3.5cm]{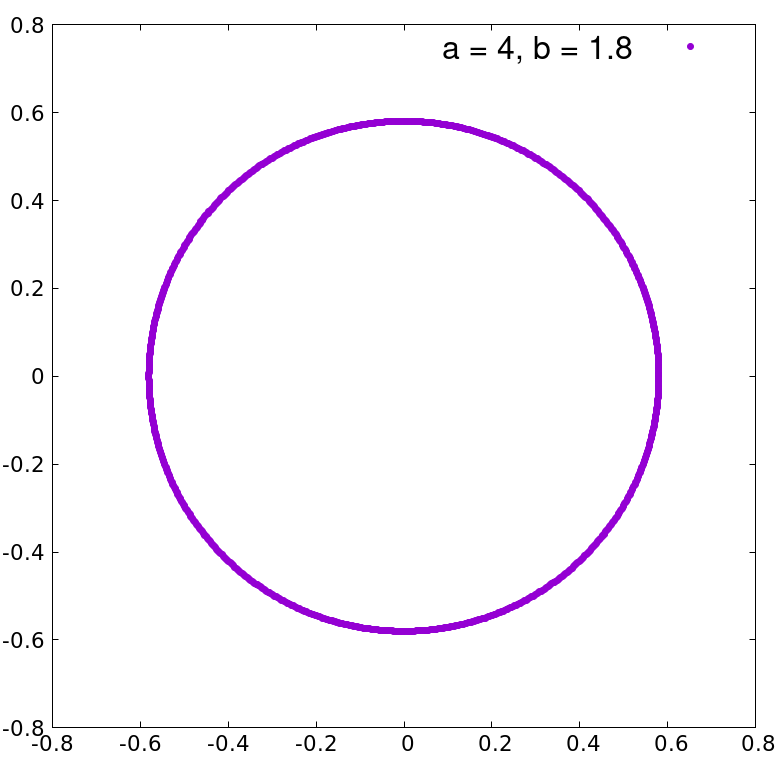}
  \caption{Minimizer for the power-law potential \eqref{eq:V-powers}
    with $a=4$ and several values of $b$, always with $10^3$
    particles.}
  \label{fig:a4}
\end{figure}

For $a=4$ the situation is similar, but now there is a range of
exponents with $b \geq 2$ where the potential $V$ is
$\mathcal{C}^2$. This is the range studied in \citet{Davies2021} and
\citet{Carrillo_2017}. In Figure \ref{fig:a4} we observe how the
increase of $b$ causes the appearance of empty space in the middle of
the particle configuration, and its further increase implies the
accumulation of all particles in a ring.

\section{Between crystallization and continuous limits}
\label{sec:between}

We will also point out that there is a range of exponents for
power-law potentials which has been left almost untouched in the
literature. For potentials $V$ which are not locally integrable (that
is, with $b \leq -d$) there are ranges of $a$ for which
crystallization (in the sense of \textbf{Q1} to \textbf{Q3})
is not expected to take place, or has not been considered. The main
unexplored range is $b \leq -d$, $a > 0$, since in the crystallization
literature and in problems motivated by statistical mechanics it is
always assumed that $V(x)$ decays to $0$ as $|x| \to +\infty$. But
there is also a range of exponents $b < a < 0$ for which
crystallization does not seem to happen, and we are left with some
kind of limiting behavior which is not yet understood. We conjecture
that in this case there is still a scaling which leads to a
well-defined limiting profile along the lines of \textbf{MO}:
\begin{description}
\item[\textbf{MS (macroscopic scaling)}] Do there exist $s > d$, a
  nontrivial measure $\mu$ on $\R^d$ and a sequence of minimizers
  $(X_N)_{N \geq 2}$ of $E_N$ such that (possibly up to a subsequence)
  \begin{equation}
    \label{eq:macro-scaling}
    \frac{1}{N} \mu_{X_N / N^{1/s}} \overset{*}{\rightharpoonup} \mu
    \qquad \text{as $N \to +\infty$?}
  \end{equation}
  Or even stronger: does this take place for \emph{all} sequences of
  minimizers, up to rigid movements of these minimizers?
\end{description}

\begin{figure}[h!]
  \centering
   \includegraphics[width=5cm]{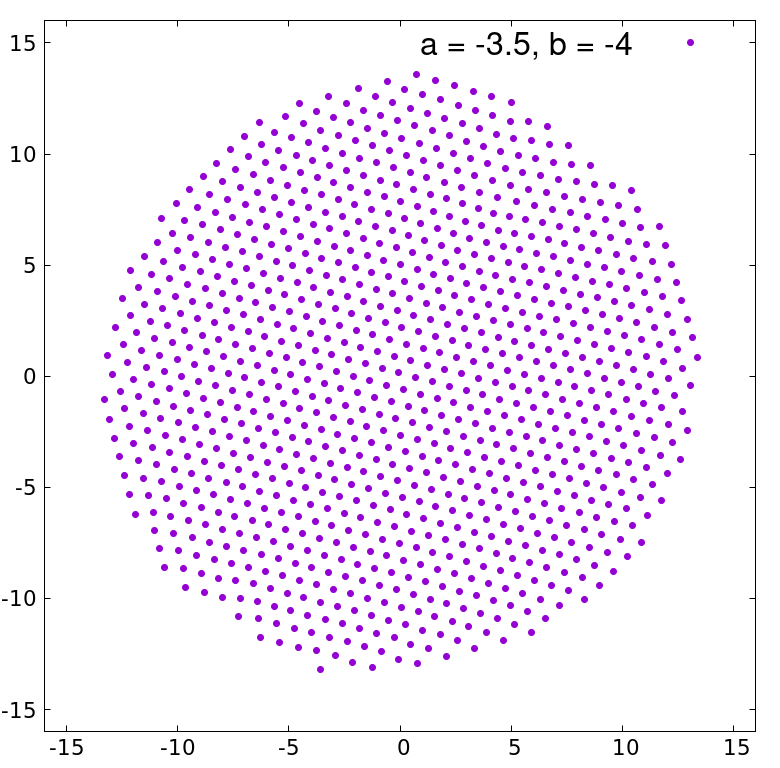}
   \includegraphics[width=5cm]{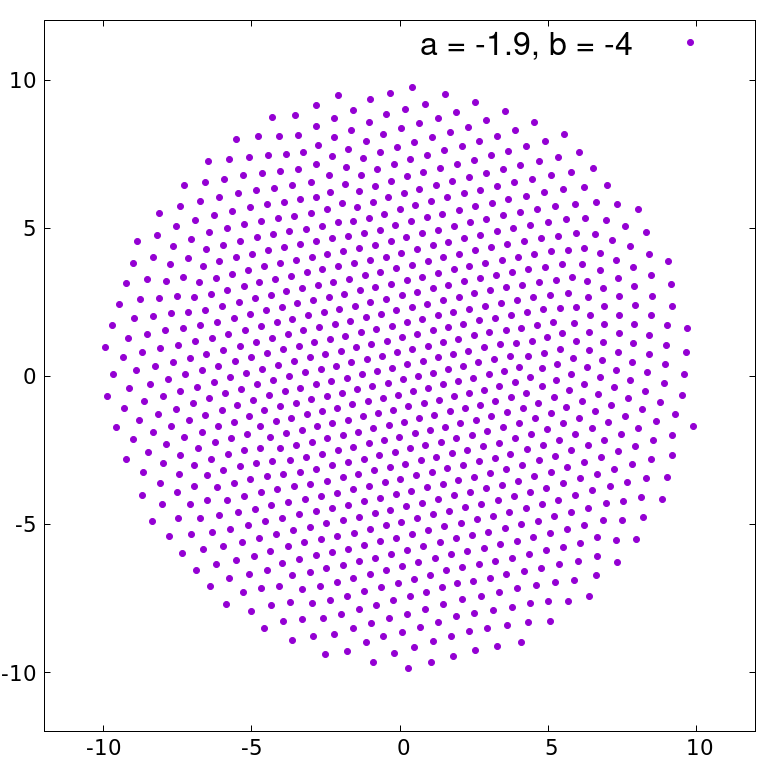}
   \includegraphics[width=5cm]{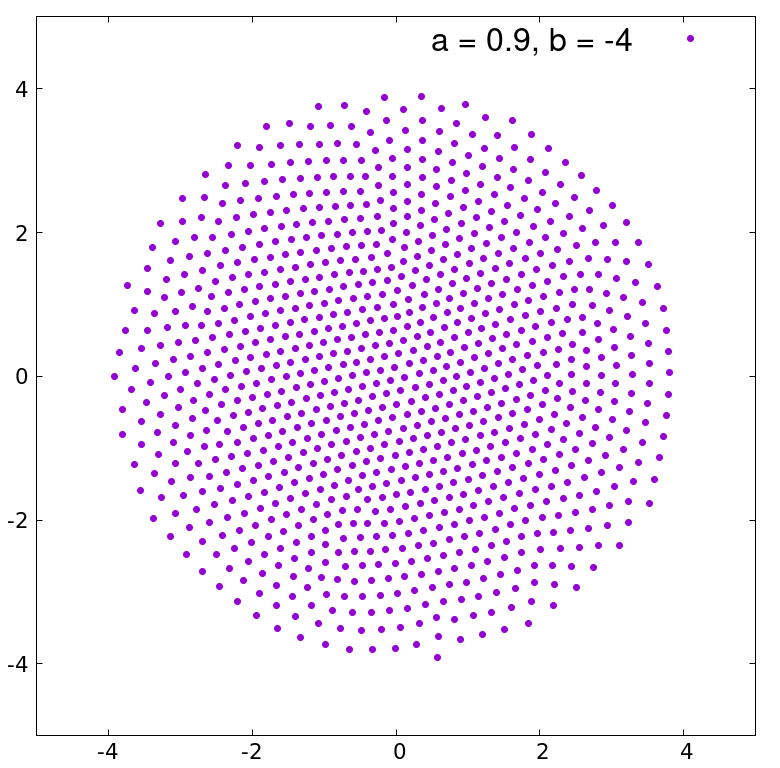}
   \includegraphics[width=5cm]{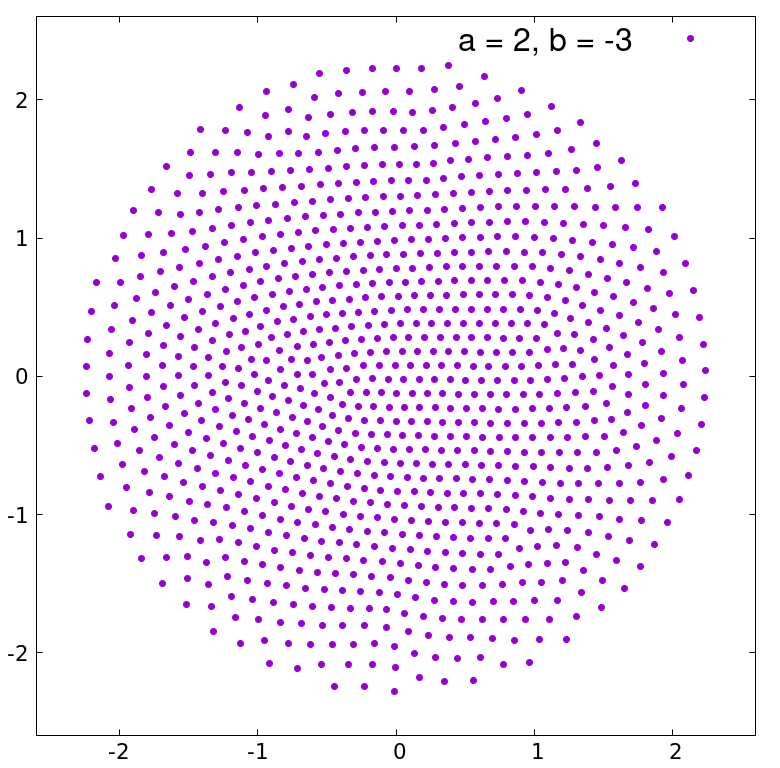}
   \includegraphics[width=5cm]{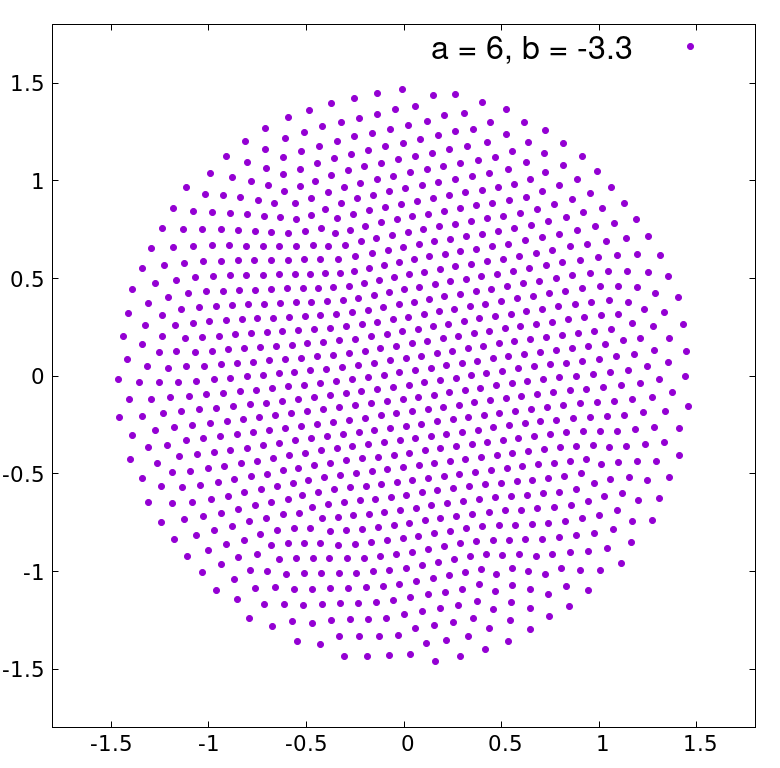}
   \includegraphics[width=5cm]{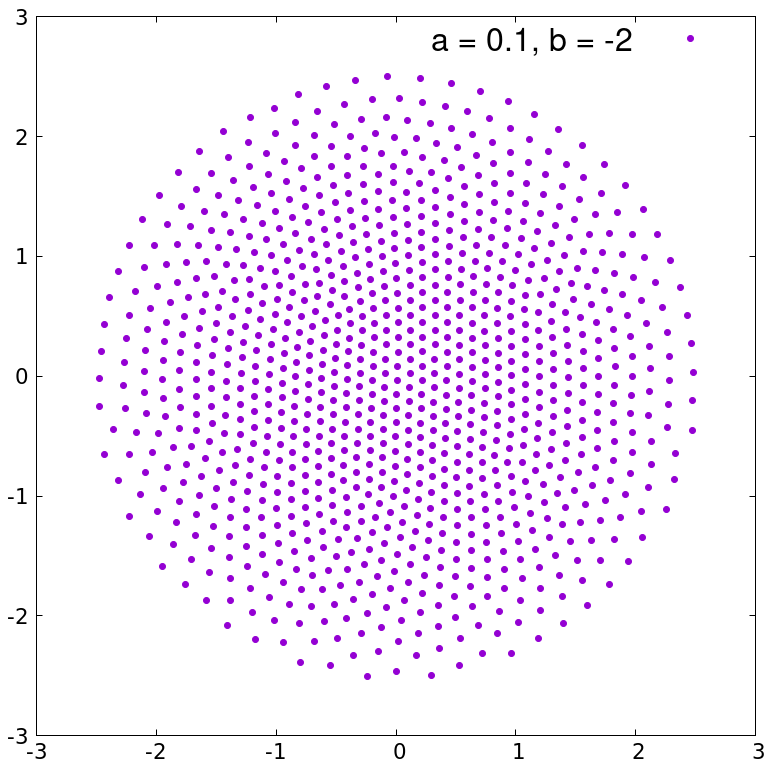}
  \caption{Some minimizers with $b \leq -2$.}
  \label{fig:bsmall}
\end{figure}

There are numerical results in 1D which seem to indicate that this is
the case \citep{Carrillo2014b}, but this has not been systematically
explored, much less the subject of any rigorous proofs. Figure
\ref{fig:bsmall} shows 2D configurations of the kind $ b<a<0 $ and
$ b<-d $, $ a>0 $. In some of them, we appreciate certain similarities
with crystallization cases (see Figure \ref{fig:Lennard-Jones}), but
we are not aware of any works exploring their behavior. This regime
may be related to scaling limits of dynamical problems with
non-integrable potentials such as \citet{Oelschlaeger1990}, but the
link is far from clear to us.

\newpage
\vspace*{6cm}
\section*{Acknowledgements}

We would like to thank David J.~Wales for several interesting discussions
and help with the software GMIN, used for the main numerical
calculations in this paper. We also thank J.~A.~Carrillo and R.~Frank
for some interesting discussions. The authors acknowledge support from
grant PID2020-117846GB-I00, the research network RED2022-134784-T, and
the María de Maeztu grant CEX2020-001105-M from the Spanish
government.

\newpage
\bibliography{canizo}

\begin{thebibliography}{74}
\providecommand{\natexlab}[1]{#1}
\providecommand{\url}[1]{\texttt{#1}}
\expandafter\ifx\csname urlstyle\endcsname\relax
  \providecommand{\doi}[1]{doi: #1}\else
  \providecommand{\doi}{doi: \begingroup \urlstyle{rm}\Url}\fi

\bibitem[Au~Yeung et~al.(2011)Au~Yeung, Friesecke, and Schmidt]{AuYeung2011}
Au~Yeung, Y., Friesecke, G., and Schmidt, B.
\newblock \href{http://dx.doi.org/10.1007/s00526-011-0427-6}{Minimizing atomic
  configurations of short range pair potentials in two dimensions:
  crystallization in the {Wulff} shape}.
\newblock \emph{Calculus of Variations and Partial Differential Equations},
  44\penalty0 (1–2):\penalty0 81--100, July 2011.

\bibitem[Balagué et~al.(2013{\natexlab{a}})Balagué, Carrillo, Laurent, and
  Raoul]{Balague2013}
Balagué, D., Carrillo, J.~A., Laurent, T., and Raoul, G.
\newblock \href{http://dx.doi.org/10.1007/s00205-013-0644-6}{{Dimensionality of
  Local Minimizers of the Interaction Energy}}.
\newblock \emph{Archive for Rational Mechanics and Analysis}, 209\penalty0
  (3):\penalty0 1055--1088, September 2013{\natexlab{a}},
  \href{http://arxiv.org/abs/1210.6795}{{\ttfamily arXiv:1210.6795}}.

\bibitem[Balagué et~al.(2013{\natexlab{b}})Balagué, Carrillo, Laurent, and
  Raoul]{Balague2013a}
Balagué, D., Carrillo, J.~A., Laurent, T., and Raoul, G.
\newblock \href{http://dx.doi.org/10.1016/j.physd.2012.10.002}{{Nonlocal
  interactions by repulsive–attractive potentials: Radial ins/stability}}.
\newblock \emph{Physica D: Nonlinear Phenomena}, 260:\penalty0 5--25, October
  2013{\natexlab{b}}, \href{http://arxiv.org/abs/1109.5258}{{\ttfamily
  arXiv:1109.5258}}.

\bibitem[Bavaud(1991)]{Bavaud1991}
Bavaud, F.
\newblock \href{http://dx.doi.org/10.1103/revmodphys.63.129}{Equilibrium
  properties of the {Vlasov} functional: The generalized
  {Poisson-Boltzmann-Emden} equation}.
\newblock \emph{Reviews of Modern Physics}, 63\penalty0 (1):\penalty0 129+,
  January 1991.

\bibitem[Ben Haj~Yedder et~al.(2003)Ben Haj~Yedder, Blanc, and
  Le~Bris]{ben2003numerical}
Ben Haj~Yedder, A., Blanc, X., and Le~Bris, C.
\newblock A numerical investigation of the 2-dimensional crystal problem.
\newblock Preprint CERMICS, available at
  http://www.ann.jussieu.fr/publications/2003/R03003.html, 2003.

\bibitem[Bertozzi et~al.(2012)Bertozzi, Laurent, and Léger]{Bertozzi2012}
Bertozzi, A.~L., Laurent, T., and Léger, F.
\newblock \href{http://dx.doi.org/10.1142/s0218202511400057}{Aggregation and
  spreading vie the {Newtonian} potential: the dynamics of patch solutions}.
\newblock \emph{Mathematical Models and Methods in the Applied Sciences},
  22\penalty0 (supp01):\penalty0 1140005+, February 2012.

\bibitem[Blanc and Lewin(2015)]{Blanc2015}
Blanc, X. and Lewin, M.
\newblock \href{http://dx.doi.org/10.4171/emss/13}{The crystallization
  conjecture: a review}.
\newblock \emph{EMS Surveys in Mathematical Sciences}, 2\penalty0 (2):\penalty0
  225--306, 2015.

\bibitem[Borodachov et~al.(2019)Borodachov, Hardin, and Saff]{Borodachov2019}
Borodachov, S.~V., Hardin, D.~P., and Saff, E.~B.
\newblock \href{http://dx.doi.org/10.1007/978-0-387-84808-2}{\emph{Discrete
  Energy on Rectifiable Sets}}.
\newblock Springer New York, 2019.

\bibitem[Bétermin(2023)]{Betermin2023}
Bétermin, L.
\newblock \href{http://dx.doi.org/10.1088/1751-8121/acc21d}{Optimality of the
  triangular lattice for {Lennard–Jones} type lattice energies: a
  computer-assisted method}.
\newblock \emph{Journal of Physics A: Mathematical and Theoretical},
  56\penalty0 (14):\penalty0 145204, March 2023.

\bibitem[Bétermin et~al.(2021)Bétermin, De~Luca, and Petrache]{B_termin_2021}
Bétermin, L., De~Luca, L., and Petrache, M.
\newblock \href{http://dx.doi.org/10.1007/s00205-021-01627-6}{Crystallization
  to the square lattice for a two-body potential}.
\newblock \emph{Archive for Rational Mechanics and Analysis}, 240\penalty0
  (2):\penalty0 987--1053, March 2021.

\bibitem[Caffarelli and V\'{a}zquez(2011{\natexlab{a}})]{Caffarelli2011}
Caffarelli, L. and V\'{a}zquez, J.~L.
\newblock \href{http://dx.doi.org/10.3934/dcds.2011.29.1393}{{Asymptotic
  behaviour of a porous medium equation with fractional diffusion}}.
\newblock \emph{Discrete and Continuous Dynamical Systems - Series A},
  29\penalty0 (4):\penalty0 1393--1404, April 2011{\natexlab{a}},
  \href{http://arxiv.org/abs/1004.1096}{{\ttfamily arXiv:1004.1096}}.

\bibitem[Caffarelli and V\'{a}zquez(2011{\natexlab{b}})]{Caffarelli2011a}
Caffarelli, L.~A. and V\'{a}zquez, J.~L.
\newblock \href{http://dx.doi.org/10.1007/s00205-011-0420-4}{{Nonlinear porous
  medium flow with fractional potential pressure}}.
\newblock \emph{Archive for Rational Mechanics and Analysis}, 202\penalty0
  (2):\penalty0 537--565, November 2011{\natexlab{b}},
  \href{http://arxiv.org/abs/1001.0410}{{\ttfamily arXiv:1001.0410}}.

\bibitem[Ca{\~{n}}izo and Patacchini(2018)]{Canizo2018c}
Ca{\~{n}}izo, J.~A. and Patacchini, F.~S.
\newblock \href{http://dx.doi.org/10.1007/s00526-017-1289-3}{Discrete
  minimisers are close to continuum minimisers for the interaction energy}.
\newblock \emph{Calculus of Variations and Partial Differential Equations},
  57\penalty0 (1), jan 2018.

\bibitem[Carrillo et~al.(2016)Carrillo, Delgadino, and Mellet]{Carrillo_2016}
Carrillo, J.~A., Delgadino, M.~G., and Mellet, A.
\newblock \href{http://dx.doi.org/10.1007/s00220-016-2598-7}{Regularity of
  local minimizers of the interaction energy via obstacle problems}.
\newblock \emph{Communications in Mathematical Physics}, 343\penalty0
  (3):\penalty0 747--781, March 2016,
  \href{http://arxiv.org/abs/1406.4040}{{\ttfamily arXiv:1406.4040}}.

\bibitem[Carrillo et~al.(2021)Carrillo, Mateu, Mora, Rondi, Scardia, and
  Verdera]{Carrillo2021a}
Carrillo, J.~A., Mateu, J., Mora, M.~G., Rondi, L., Scardia, L., and Verdera,
  J.
\newblock \href{http://dx.doi.org/10.1007/s00526-021-01928-4}{The equilibrium
  measure for an anisotropic nonlocal energy}.
\newblock \emph{Calculus of Variations and Partial Differential Equations},
  60\penalty0 (3), 2021.

\bibitem[Carrillo et~al.(2017)Carrillo, Figalli, and Patacchini]{Carrillo_2017}
Carrillo, J., Figalli, A., and Patacchini, F.
\newblock \href{http://dx.doi.org/10.1016/j.anihpc.2016.10.004}{Geometry of
  minimizers for the interaction energy with mildly repulsive potentials}.
\newblock \emph{Annales de l’Institut Henri Poincaré C, Analyse non
  linéaire}, 34\penalty0 (5):\penalty0 1299--1308, October 2017.

\bibitem[Carrillo and Shu(2023)]{Carrillo2023a}
Carrillo, J.~A. and Shu, R.
\newblock \href{http://dx.doi.org/10.1002/cpa.22162}{Global minimizers of a
  large class of anisotropic attractive‐repulsive interaction energies in
  {2D}}.
\newblock \emph{Communications on Pure and Applied Mathematics}, 77\penalty0
  (2):\penalty0 1353--1404, September 2023.

\bibitem[Carrillo and Huang(2017)]{A_Carrillo_2017}
Carrillo, J.~A. and Huang, Y.
\newblock \href{http://dx.doi.org/10.3934/krm.2017007}{Explicit equilibrium
  solutions for the aggregation equation with power-law potentials}.
\newblock \emph{Kinetic and Related Models}, 10\penalty0 (1):\penalty0
  171--192, 2017.

\bibitem[Carrillo and Shu(2021)]{Carrillo_2022}
Carrillo, J.~A. and Shu, R.
\newblock \href{http://dx.doi.org/10.1007/s00526-022-02368-4}{From radial
  symmetry to fractal behavior of aggregation equilibria for
  repulsive–attractive potentials}.
\newblock \emph{Calculus of Variations and Partial Differential Equations},
  62\penalty0 (1), November 2021.

\bibitem[Carrillo and Shu(2022)]{Carrillo2022}
Carrillo, J.~A. and Shu, R.
\newblock \href{http://dx.doi.org/10.1515/acv-2022-0059}{Minimizers of {3D}
  anisotropic interaction energies}.
\newblock \emph{Advances in Calculus of Variations}, November 2022.

\bibitem[Carrillo et~al.(2014)Carrillo, Chipot, and Huang]{Carrillo2014b}
Carrillo, J.~A., Chipot, M., and Huang, Y.
\newblock \href{http://dx.doi.org/10.1098/rsta.2013.0399}{On global minimizers
  of repulsive–attractive power-law interaction energies}.
\newblock \emph{Philosophical Transactions of the Royal Society A:
  Mathematical, Physical and Engineering Sciences}, 372\penalty0
  (2028):\penalty0 20130399, 2014.

\bibitem[Cañizo et~al.(2015)Cañizo, Carrillo, and Patacchini]{Canizo2015b}
Cañizo, J.~A., Carrillo, J.~A., and Patacchini, F.~S.
\newblock \href{http://dx.doi.org/10.1007/s00205-015-0852-3}{Existence of
  compactly supported global minimisers for the interaction energy}.
\newblock \emph{Archive for Rational Mechanics and Analysis}, 217\penalty0
  (3):\penalty0 1197--1217, March 2015,
  \href{http://arxiv.org/abs/1405.5428}{{\ttfamily arXiv:1405.5428}}.

\bibitem[Cohn et~al.(2017)Cohn, Kumar, Miller, Radchenko, and
  Viazovska]{Cohn2017}
Cohn, H., Kumar, A., Miller, S., Radchenko, D., and Viazovska, M.
\newblock \href{http://dx.doi.org/10.4007/annals.2017.185.3.8}{The sphere
  packing problem in dimension $24$}.
\newblock \emph{Annals of Mathematics}, 185\penalty0 (3), 2017.

\bibitem[Cohn et~al.(2022)Cohn, Kumar, Miller, Radchenko, and
  Viazovska]{Cohn2022}
Cohn, H., Kumar, A., Miller, S., Radchenko, D., and Viazovska, M.
\newblock \href{http://dx.doi.org/10.4007/annals.2022.196.3.3}{Universal
  optimality of the {$E_8$} and {Leech} lattices and interpolation formulas}.
\newblock \emph{Annals of Mathematics}, 196\penalty0 (3), 2022.

\bibitem[Davies et~al.(2021)Davies, Lim, and McCann]{Davies2021}
Davies, C., Lim, T., and McCann, R.~J.
\newblock \href{http://dx.doi.org/10.1137/21m1455309}{Classifying minimum
  energy states for interacting particles: Spherical shells}.
\newblock \emph{SIAM Journal on Applied Mathematics}, 82\penalty0 (4):\penalty0
  1520--1536, August 2021.

\bibitem[Davies et~al.(2022)Davies, Lim, and McCann]{Davies_2022}
Davies, C., Lim, T., and McCann, R.~J.
\newblock \href{http://dx.doi.org/10.1007/s00220-022-04564-x}{Classifying
  minimum energy states for interacting particles: Regular simplices}.
\newblock \emph{Communications in Mathematical Physics}, 399\penalty0
  (2):\penalty0 577--598, November 2022.

\bibitem[Deaven et~al.(1996)Deaven, Tit, Morris, and Ho]{Daven1996}
Deaven, D.~M., Tit, N., Morris, J.~R., and Ho, K.~M.
\newblock \href{http://dx.doi.org/10.1016/0009-2614(96)00406-x}{Structural
  optimization of {Lennard-Jones} clusters by a genetic algorithm}.
\newblock \emph{Chemical Physics Letters}, 256\penalty0 (1-2):\penalty0
  195--200, 1996.

\bibitem[Dittner and Hartke(2017)]{Dittner2017}
Dittner, M. and Hartke, B.
\newblock \href{http://dx.doi.org/10.1016/j.comptc.2016.09.032}{Conquering the
  hard cases of {Lennard-Jones} clusters with simple recipes}.
\newblock \emph{Computational and Theoretical Chemistry}, 1107:\penalty0 7--13,
  May 2017.

\bibitem[D'Orsogna et~al.(2006)D'Orsogna, Chuang, Bertozzi, and
  Chayes]{DOrsogna2006}
D'Orsogna, M.~R., Chuang, Y.~L., Bertozzi, A.~L., and Chayes, L.~S.
\newblock \href{http://dx.doi.org/10.1103/physrevlett.96.104302}{Self-propelled
  particles with soft-core interactions: Patterns, stability, and collapse}.
\newblock \emph{Physical Review Letters}, 96\penalty0 (10):\penalty0 104302+,
  March 2006.

\bibitem[E and Li(2008)]{Weinan2008}
E, W. and Li, D.
\newblock \href{http://dx.doi.org/10.1007/s00220-008-0586-2}{On the
  crystallization of {2D} hexagonal lattices}.
\newblock \emph{Communications in Mathematical Physics}, 286\penalty0
  (3):\penalty0 1099--1140, July 2008.

\bibitem[Fellner and Raoul(2010{\natexlab{a}})]{Fellner2010}
Fellner, K. and Raoul, G.
\newblock \href{http://dx.doi.org/10.1016/j.mcm.2010.03.021}{{Stability of
  stationary states of non-local equations with singular interaction
  potentials}}.
\newblock \emph{Mathematical and Computer Modelling}, March 2010{\natexlab{a}}.

\bibitem[Fellner and Raoul(2010{\natexlab{b}})]{Fellner2010a}
Fellner, K. and Raoul, G.
\newblock \href{http://dx.doi.org/10.1142/s0218202510004921}{Stable stationary
  states of non-local interaction equations}.
\newblock \emph{Mathematical Models and Methods in Applied Sciences},
  20\penalty0 (12):\penalty0 2267--2291, December 2010{\natexlab{b}}.

\bibitem[Frank(2022)]{Frank2022}
Frank, R.~L.
\newblock \href{http://dx.doi.org/10.1016/j.na.2021.112691}{Minimizers for a
  one-dimensional interaction energy}.
\newblock \emph{Nonlinear Analysis}, 216:\penalty0 112691, 2022.

\bibitem[Frank(2023)]{Frank_2023}
Frank, R.~L.
\newblock \href{http://dx.doi.org/10.4171/8ecm/06}{\emph{Some minimization
  problems for mean field models with competing forces}}, pages 277--294.
\newblock EMS Press, July 2023.
\newblock ISBN 9783985475513.

\bibitem[Frank and Matzke(2023)]{frank2023minimizers}
Frank, R.~L. and Matzke, R.~W.
\newblock Minimizers for an aggregation model with attractive-repulsive
  interaction.
\newblock 2023, \href{http://arxiv.org/abs/2307.13769}{{\ttfamily
  arXiv:2307.13769}}.

\bibitem[Gardner and Radin(1979)]{Gardner1979}
Gardner, C.~S. and Radin, C.
\newblock \href{http://dx.doi.org/10.1007/bf01009521}{The infinite-volume
  ground state of the {Lennard-Jones} potential}.
\newblock \emph{Journal of Statistical Physics}, 20\penalty0 (6):\penalty0
  719--724, June 1979.

\bibitem[Hales(2005)]{Hales2005}
Hales, T.
\newblock \href{http://dx.doi.org/10.4007/annals.2005.162.1065}{A proof of the
  {Kepler} conjecture}.
\newblock \emph{Annals of Mathematics}, 162\penalty0 (3):\penalty0 1065--1185,
  2005.

\bibitem[Hardin and Saff(2005)]{Hardin2005}
Hardin, D.~P. and Saff, E.~B.
\newblock \href{http://dx.doi.org/10.1016/j.aim.2004.05.006}{Minimal riesz
  energy point configurations for rectifiable $d$-dimensional manifolds}.
\newblock \emph{Advances in Mathematics}, 193\penalty0 (1):\penalty0 174--204,
  2005.

\bibitem[Hayes(2012)]{Hayes2012}
Hayes, B.
\newblock \href{http://dx.doi.org/10.1511/2012.99.442}{The science of sticky
  spheres}.
\newblock \emph{American Scientist}, 100\penalty0 (6):\penalty0 442, 2012.

\bibitem[Heitmann and Radin(1980)]{Heitmann1980}
Heitmann, R.~C. and Radin, C.
\newblock \href{http://dx.doi.org/10.1007/bf01014644}{The ground state for
  sticky disks}.
\newblock \emph{Journal of Statistical Physics}, 22\penalty0 (3):\penalty0
  281--287, 1980.

\bibitem[Hoare(1979)]{Hoare1979}
Hoare, M.~R.
\newblock Structure and dynamics of simple microclusters.
\newblock In Prigogine, I. and Rice, S.~A., editors, \emph{Advances in Chemical
  Physics}, volume~40, pages 49--135. John Wiley \& Sons, 1979.

\bibitem[Hoare and Pal(1971{\natexlab{a}})]{Hoare1971}
Hoare, M.~R. and Pal, P.
\newblock \href{http://dx.doi.org/10.1080/00018737100101231}{Physical cluster
  mechanics: Statics and energy surfaces for monatomic systems}.
\newblock \emph{Advances in Chemical Physics}, 20:\penalty0 161--196,
  1971{\natexlab{a}}.

\bibitem[Hoare and Pal(1971{\natexlab{b}})]{Hoare1971a}
Hoare, M.~R. and Pal, P.
\newblock \href{http://dx.doi.org/10.1038/physci230005a0}{Statics and stability
  of small cluster nuclei}.
\newblock \emph{Nature (Physical Sciences)}, 230:\penalty0 5--8,
  1971{\natexlab{b}}.

\bibitem[Hoare and Pal(1972)]{HOARE1972}
Hoare, M.~R. and Pal, P.
\newblock \href{http://dx.doi.org/10.1038/physci236035a0}{Geometry and
  stability of ``spherical'' f.c.c. microcrystallites}.
\newblock \emph{Nature (Physical Sciences)}, 236:\penalty0 35--37, 1972.

\bibitem[Hoy et~al.(2012)Hoy, Harwayne-Gidansky, and O’Hern]{Hoy2012}
Hoy, R.~S., Harwayne-Gidansky, J., and O’Hern, C.~S.
\newblock \href{http://dx.doi.org/10.1103/physreve.85.051403}{Structure of
  finite sphere packings via exact enumeration: Implications for colloidal
  crystal nucleation}.
\newblock \emph{Physical Review E}, 85\penalty0 (5):\penalty0 051403, 2012.

\bibitem[Kirkpatrick et~al.(1983)Kirkpatrick, Gelatt, and
  Vecchi]{Kirkpatrick1983}
Kirkpatrick, S., Gelatt, C.~D., and Vecchi, M.~P.
\newblock \href{http://dx.doi.org/10.1126/science.220.4598.671}{Optimization by
  simulated annealing}.
\newblock \emph{Science}, 220\penalty0 (4598):\penalty0 671--680, 1983.

\bibitem[Kolokolnikov et~al.(2011)Kolokolnikov, Sun, Uminsky, and
  Bertozzi]{Kolokolnikov2011}
Kolokolnikov, T., Sun, H., Uminsky, D., and Bertozzi, A.~L.
\newblock \href{http://dx.doi.org/10.1103/physreve.84.015203}{{Stability of
  ring patterns arising from two-dimensional particle interactions}}.
\newblock \emph{Physical Review E}, 84:\penalty0 015203+, July 2011.

\bibitem[Li and Scheraga(1987)]{Li1987}
Li, Z. and Scheraga, H.~A.
\newblock \href{http://dx.doi.org/10.1073/pnas.84.19.6611}{{Monte
  Carlo}-minimization approach to the multiple-minima problem in protein
  folding.}
\newblock \emph{Proceedings of the National Academy of Sciences}, 84\penalty0
  (19):\penalty0 6611--6615, 1987.

\bibitem[Liu and Nocedal(1989)]{liu1989limited}
Liu, D.~C. and Nocedal, J.
\newblock On the limited memory {BFGS} method for large scale optimization.
\newblock \emph{Mathematical programming}, 45\penalty0 (1):\penalty0 503--528,
  1989.

\bibitem[Luca and Friesecke(2018)]{Luca2018}
Luca, L.~D. and Friesecke, G.
\newblock \href{http://dx.doi.org/10.1007/s00332-017-9401-6}{Crystallization in
  two dimensions and a discrete {Gauss–Bonnet} theorem}.
\newblock \emph{Journal of Nonlinear Science}, 28\penalty0 (1):\penalty0
  69--90, 2018.

\bibitem[Meng et~al.(2010)Meng, Arkus, Brenner, and Manoharan]{Meng2010}
Meng, G., Arkus, N., Brenner, M.~P., and Manoharan, V.~N.
\newblock \href{http://dx.doi.org/10.1126/science.1181263}{The free-energy
  landscape of clusters of attractive hard spheres}.
\newblock \emph{Science}, 327\penalty0 (5965):\penalty0 560--563, 2010.

\bibitem[Niesse and Mayne(1996)]{Niesse1996}
Niesse, J.~A. and Mayne, H.~R.
\newblock \href{http://dx.doi.org/10.1016/0009-2614(96)01000-7}{Minimization of
  small silicon clusters using the space-fixed modified genetic algorithm
  method}.
\newblock \emph{Chemical Physics Letters}, 261\penalty0 (4-5):\penalty0
  576--582, 1996.

\bibitem[Northby(1987)]{Northby1987}
Northby, J.~A.
\newblock \href{http://dx.doi.org/10.1063/1.453492}{Structure and binding of
  {Lennard-Jones} clusters: {$13 \leq N \leq 147$}}.
\newblock \emph{Journal of Chemical Physics}, 87\penalty0 (10):\penalty0
  6166--6177, 1987.

\bibitem[Oelschläger(1990)]{Oelschlaeger1990}
Oelschläger, K.
\newblock \href{http://dx.doi.org/10.1016/0022-0396(90)90101-t}{Large systems
  of interacting particles and the porous medium equation}.
\newblock \emph{Journal of Differential Equations}, 88\penalty0 (2):\penalty0
  294--346, December 1990.

\bibitem[Petrache and Serfaty(2017)]{Petrache2017}
Petrache, M. and Serfaty, S.
\newblock \href{http://dx.doi.org/10.1017/s1474748015000201}{Next order
  asymptotics and renormalized energy for {Riesz} interactions}.
\newblock \emph{Journal of the Institute of Mathematics of Jussieu},
  16\penalty0 (3):\penalty0 501--569, 2017.

\bibitem[Petrache and Serfaty(2020)]{Petrache}
Petrache, M. and Serfaty, S.
\newblock \href{http://dx.doi.org/10.1090/proc/15003}{Crystallization for
  {Coulomb} and {Riesz} interactions as a consequence of the {Cohn-Kumar}
  conjecture}.
\newblock \emph{Proceedings of the American Mathematical Society}, 148\penalty0
  (7):\penalty0 3047--3057, 2020.

\bibitem[Radin(1981)]{Radin1981}
Radin, C.
\newblock \href{http://dx.doi.org/10.1007/bf01013177}{{The ground state for
  soft disks}}.
\newblock \emph{Journal of Statistical Physics}, 26\penalty0 (2):\penalty0
  365--373, October 1981.

\bibitem[Radin(1987)]{Radin1987}
Radin, C.
\newblock \href{http://dx.doi.org/10.1142/S0217979287001675}{{Low temperature
  and the origin of crystalline symmetry}}.
\newblock \emph{International Journal of Modern Physics B}, 1\penalty0
  (5/6):\penalty0 1157--1191, October 1987.

\bibitem[Raoul(2012)]{Raoul2012}
Raoul, G.
\newblock \href{http://dx.doi.org/10.57262/die/1356012673}{Nonlocal interaction
  equations: Stationary states and stability analysis}.
\newblock \emph{Differential and Integral Equations}, 25\penalty0 (5/6), 2012.

\bibitem[Ruelle(1969)]{Ruelle1969}
Ruelle, D.
\newblock \href{http://www.worldcat.org/isbn/9780805383607}{\emph{{Statistical
  mechanics: rigorous results}}}.
\newblock W.A. Benjamin, 1969.
\newblock ISBN 9780805383607.

\bibitem[Simione et~al.(2015)Simione, Slep\v{c}ev, and Topaloglu]{Simione2015}
Simione, R., Slep\v{c}ev, D., and Topaloglu, I.
\newblock \href{http://dx.doi.org/10.1007/s10955-015-1215-z}{{Existence of
  Ground States of Nonlocal-Interaction Energies}}.
\newblock \emph{Journal of Statistical Physics}, 159\penalty0 (4):\penalty0
  972--986, February 2015, \href{http://arxiv.org/abs/1405.5146}{{\ttfamily
  arXiv:1405.5146}}.

\bibitem[Theil(2006)]{Theil2006}
Theil, F.
\newblock \href{http://dx.doi.org/10.1007/s00220-005-1458-7}{{A Proof of
  Crystallization in Two Dimensions}}.
\newblock \emph{Communications in Mathematical Physics}, 262\penalty0
  (1):\penalty0 209--236, February 2006.

\bibitem[Trombach et~al.(2018)Trombach, Hoy, Wales, and
  Schwerdtfeger]{Trombach}
Trombach, L., Hoy, R.~S., Wales, D.~J., and Schwerdtfeger, P.
\newblock \href{http://dx.doi.org/10.1103/physreve.97.043309}{From
  sticky-hard-sphere to {Lennard-Jones-type} clusters}.
\newblock \emph{Physical Review E}, 97\penalty0 (4):\penalty0 043309, 2018.

\bibitem[Ventevogel(1978{\natexlab{a}})]{Ventevogel1978}
Ventevogel, W.~J.
\newblock \href{http://dx.doi.org/10.1016/0375-9601(78)90685-0}{Why do crystals
  exist?}
\newblock \emph{Physics Letters A}, 64\penalty0 (5):\penalty0 463--464,
  1978{\natexlab{a}}.

\bibitem[Ventevogel(1978{\natexlab{b}})]{Ventevogel1978a}
Ventevogel, W.~J.
\newblock \href{http://dx.doi.org/10.1016/0378-4371(78)90136-x}{On the
  configuration of a one-dimensional system of interacting particles with
  minimum potential energy per particle}.
\newblock \emph{Physica A: Statistical Mechanics and its Applications},
  92\penalty0 (3-4):\penalty0 343--361, 1978{\natexlab{b}}.

\bibitem[Ventevogel and Nijboer(1979{\natexlab{a}})]{Ventevogel1979}
Ventevogel, W.~J. and Nijboer, B. R.~A.
\newblock \href{http://dx.doi.org/10.1016/0378-4371(79)90072-4}{On the
  configuration of systems of interacting particles with minimum potential
  energy per particle}.
\newblock \emph{Physica A: Statistical Mechanics and its Applications},
  99\penalty0 (3):\penalty0 569--580, 1979{\natexlab{a}}.

\bibitem[Ventevogel and Nijboer(1979{\natexlab{b}})]{Ventevogel1979a}
Ventevogel, W.~J. and Nijboer, B. R.~A.
\newblock \href{http://dx.doi.org/10.1016/0378-4371(79)90178-x}{On the
  configuration of systems of interacting particle with minimum potential
  energy per particle}.
\newblock \emph{Physica A: Statistical Mechanics and its Applications},
  98\penalty0 (1-2):\penalty0 274--288, 1979{\natexlab{b}}.

\bibitem[Viazovska(2017)]{Viazovska2017}
Viazovska, M.
\newblock \href{http://dx.doi.org/10.4007/annals.2017.185.3.7}{The sphere
  packing problem in dimension $8$}.
\newblock \emph{Annals of Mathematics}, 185\penalty0 (3), 2017.

\bibitem[Wales et~al.()Wales, Doye, Dullweber, Hodges, Naumkin, Calvo,
  Hernández-Rojas, and Middleton]{CCDonline}
Wales, D.~J., Doye, J. P.~K., Dullweber, A., Hodges, M.~P., Naumkin, F.~Y.,
  Calvo, F., Hernández-Rojas, J., and Middleton, T.~F.
\newblock \href{http://www-wales.ch.cam.ac.uk/CCD.html}{The {Cambridge}
  {Cluster} {Database}}.
\newblock Online.
\newblock URL http://www-wales.ch.cam.ac.uk/CCD.html.

\bibitem[Wales()]{WalesGMIN}
Wales, D.~J.
\newblock \href{https://www-wales.ch.cam.ac.uk/software.html}{{GMIN}
  optimisation software}.
\newblock Available at https://www-wales.ch.cam.ac.uk/GMIN/.

\bibitem[Wales and Doye(1997)]{Wales1997}
Wales, D.~J. and Doye, J. P.~K.
\newblock \href{http://dx.doi.org/10.1021/jp970984n}{Global optimization by
  basin-hopping and the lowest energy structures of {Lennard-Jones} clusters
  containing up to 110 atoms}.
\newblock \emph{The Journal of Physical Chemistry A}, 101\penalty0
  (28):\penalty0 5111--5116, 1997.

\bibitem[Wille(1987{\natexlab{a}})]{Wille1987}
Wille, L.~T.
\newblock \href{http://dx.doi.org/10.1038/325374c0}{Searching potential energy
  surfaces by simulated annealing}.
\newblock \emph{Nature}, 325:\penalty0 374--374, 1987{\natexlab{a}}.

\bibitem[Wille(1987{\natexlab{b}})]{Wille1987a}
Wille, L.~T.
\newblock \href{http://dx.doi.org/10.1016/0009-2614(87)87091-4}{Minimum-energy
  configurations of atomic clusters: new results obtained by simulated
  annealing}.
\newblock \emph{Chemical Physics Letters}, 133\penalty0 (5):\penalty0 405--410,
  January 1987{\natexlab{b}}.

\bibitem[Zhang(2010)]{Zhang2010}
Zhang, J.
\newblock A brief review on results and computational algorithms for minimizing
  the {Lennard-Jones} potential.
\newblock December 2010, \href{http://arxiv.org/abs/1101.0039}{{\ttfamily
  arXiv:1101.0039 [physics.comp-ph]}}.

\end{thebibliography}

\end{document}